\documentclass[12pt]{amsart}
\usepackage[all]{xy}
\usepackage{amssymb}
\usepackage{amsopn}
\usepackage{amstext}
\usepackage{amsmath}

\newtheorem{defn}[equation]{Definition}
\newtheorem{thm}[equation]{Theorem}
\newtheorem{thmSN}{Theorem}
\newtheorem{prop}[equation]{Proposition}{\bf}{\it}
\newtheorem{rem}[equation]{Remark}{\it}{\rm}
\newtheorem{exam}[equation]{Example}{\it}{\rm}

\numberwithin{equation}{section}

\renewcommand\sf{\operatorname{Spf}}
\newcommand\grv{{\operatorname{Gr}}(V)}
\newcommand\gr{\operatorname{Gr}}

\newcommand\nm{\operatorname{Nm}}

\newcommand\Det{\operatorname{Det}}

\newcommand\Detd{\operatorname{Det}^\ast}

\renewcommand\ker{\operatorname{Ker}}
\DeclareMathOperator{\res}{{Res}}

\newcommand\limil[1]{\underset{#1}\varinjlim\,}
\newcommand\limp{\varprojlim}

\newcommand\aut{\operatorname{Aut}}

\newcommand\End{\operatorname{End}}
\newcommand\tr{\operatorname{Tr}}

\newcommand\rank{\operatorname{rank}}

\newcommand\pic{\operatorname{Pic}}

\newcommand\kr{\operatorname{Kr}}

\renewcommand\o{{\mathcal O}}
\renewcommand\j{{\mathcal J}}

\renewcommand\L{{\mathcal L}}

\newcommand\M{{\mathcal M}}

\newcommand\Z{{\mathbb Z}}

\newcommand\C{{\mathbb C}}

\newcommand\w{\widehat}
\newcommand\wtilde{\widetilde}
\newcommand\iso{\overset{\sim}\to}

\begin{document}

\title[The Schottky problem for curves with automorphisms]{A solution of the Schottky-Type problem \\ for curves with
automorphisms}
\author[E. G\'omez \and J. M. Mu\~noz \and F. J. Plaza \and S. Recillas \and R. E. Rodr\'{\i}guez]
{Esteban G\'omez Gonz\'alez${\;}^{1}$ \\ Jos\'e M. Mu\~noz
Porras${\;}^{1}$ \\  Francisco J. Plaza Mart\'{\i}n${\;}^{1}$ \\ Sev\'{\i}n
Recillas${\;}^{2}$ \\ Rub\'{\i} E. Rodr\'{\i}guez${\;}^{3}$ }

\thanks{This work is partially supported by the research contracts
BFM2003-00078 of DGI and  SA071/04 of JCyL, CONACYT grant 40033-F
and Fondecyt grant 1030595. The third author is also supported by
MCYT ``Ram\'on y Cajal'' program.%
\\
\indent {\it Keywords:} Schottky problem, curves with
automorphisms, theta
functions, KP and KdV hierarchies, completely integrable systems.} %

\address{${\!}^{1}$ Departamento de Matem\'aticas, Universidad de Salamanca,
Plaza de la Merced 1-4, 37008 Salamanca. Spain.
          \\
${\!}^{2}$ Instituto de Matem\'aticas, UNAM Campus Morelia,
Morelia Mich., 58089, M\'{e}xico, and  CIMAT, Callej\'{o}n de Jalisco s/n,
Valenciana Gto., 36000, M\'{e}xico.
        \\
${\!}^{3}$ Departamento de Matem\'aticas, Pontificia Universidad
Cat\'olica de Chile, Casilla 306-22, Santiago. Chile. }

\email{esteban@usal.es} \email{jmp@usal.es} \email{fplaza@usal.es}
\email{sevin@matmor.unam.mx} \email{sevin@cimat.mx}
\email{rubi@mat.puc.cl}

\subjclass[2000]{Primary: 14H42, 14H37; Secondary: 37K10, 14K25}



%
%
\begin{abstract}
In this paper, an explicit hierarchy of differential equations for
the $\tau$-functions defining the moduli space of curves with
automorphisms as a subscheme of the Sato Grassmannian is obtained.
The Schottky problem for Riemann surfaces with automorphisms
consists of characterizing those p.p.a.v. that are Jacobian
varieties of a curve with a non-trivial automorphism. A
characterization in terms of hierarchies of p.d.e. for theta
functions is also given.
\end{abstract}


\maketitle



\section{Introduction}

The objective of the Schottky problem is to characterize the
principally polarized abelian varieties (p.p.a.v.) which are
Jacobians of smooth algebraic curves. This problem was solved, in
the framework of the theory of KP equations, by Shiota
(\cite{Shiota}). The analogous problem for Prym varieties was
studied and partially solved by Shiota (\cite{Sh2}) in terms of
the BKP hierarchy.

The Schottky problem for Pryms is also related to the
characterization of Jacobians of algebraic curves which admit a
non-trivial involution. The moduli space of curves with
non-trivial automorphisms was studied in the previous paper
(\cite{GMP}), and the points of the Sato Grassmannian defined by
those curves were characterized. Moreover, an explicit set of
algebraic equations defining the moduli space of curves with
automorphisms as a subscheme of the Sato Grassmannian were
obtained.

In \S3 of the present paper, these equations are given as an
explicit hierarchy of differential equations for the
$\tau$-functions (Theorems~\ref{thm:equationsMram}
and~\ref{thm:equationsMNram}). The first non-trivial equations are
studied in detail and related to the KP and KdV hierarchies.

These results do not solve the Schottky problem for Riemann
surfaces with automorphisms. For solving this problem one must
characterize the conditions that a theta function of a p.p.a.v.
should satisfy in order to be a jacobian theta function of a curve
with a non-trivial automorphism. Such a characterization is
obtained in this paper, as follows.

Let $X_\Omega$ be an irreducible p.p.a.v. of dimension $g$. For a
natural number $r$, the $r$-th Baker-Akhiezer functions are
defined from the theta function of $X_\Omega$ and some auxiliary
data; namely, an $r$-tuple of $g \times \infty$-matrices
$(A^{(1)},\dots, A^{(r)})$ of rank $g$ and $r$ symmetric quadratic
forms $Q^{(1)},\dots, Q^{(r)}$ (see subsection~\ref{sec:taujac}
for precise definitions).

\begin{thmSN}[Characterization]
Let $X_\Omega$ be an irreducible p.p.a.v. of dimension $g>1$.

Then the following conditions are equivalent.
\begin{enumerate}
    \item There exists a projective irreducible smooth genus $g$ curve $C$ with a non-trivial
    automorphism $\sigma_C:C\to C$ such that $X_\Omega$ is
    isomorphic as a p.p.a.v. to the Jacobian of $C$.

    \item There exist a prime number $p$, $p$ matrices $A^{(1)},\dots, A^{(p)}$ ($A^{(j)}$
    being a $g \times \infty$-matrix of rank $g$) and $p$ symmetric quadratic
    forms $Q^{(1)},\dots, Q^{(p)}$, such that
        \begin{enumerate}
        \item[a)] for some $\xi_0\in \C^g$ the corresponding
        BA-functions satisfy the $(1,\overset{p}{\ldots}, 1)$-KP
        hierarchy
            $$
            \res \left( \sum_{j=1}^p z^{-\delta_{ju} - \delta_{jv}} \,
            \psi_{u,\xi_0}^{(j)} (z,t) \psi_{v,\xi_0}^{*(j)} (z,s) \right) \, dz
            \,=\,0 \ , \text{ and }
            $$
        \item[b)] there exists $\xi_1\in \C^g$ (depending on $\xi_0$) such that
            $$
            \res \left( \sum_{j=1}^p z^{-\delta_{ju} - \delta_{jv}} \,
            \psi_{v+1,\xi_0}^{(j+1)} (z,\sigma^*(t))
            \psi_{u,\xi_1}^{*(j)} (z,s)\right) \, dz
            \,=\,0
            $$
            where $\sigma^*(t):=
            (t^{(p)}, t^{(1)}, t^{(2)}, \ldots , t^{(p-1)})$.
        \end{enumerate}
\end{enumerate}
\end{thmSN}

This theorem can be understood as a translation into equations of
the characterization theorems of~\S5, where they are stated in
terms of orbits of finite dimension as in the approach of Mulase
(\cite{Mulase}) and Shiota (\cite{Shiota}). The idea of our proof
is similar to the first part of the paper of Shiota; in fact, our
main ingredient is a generalization of Theorem~6 of Shiota to the
case of the $(1,\overset{r}{\ldots}, 1)$-KP hierarchy (see
Theorem~\ref{thm:S-Mgen}).

Standard arguments allow us to reinterpret the above identities as
a hierarchy of PDE for the associated tau functions. The final
step of our program should consist of proving that the ``first''
non-trivial equations of the hierarchy suffice to characterize the
theta functions of Riemann surfaces with non-trivial
automorphisms. This problem has an analytical counterpart which is
the analogue of Shiota's proof of Novikov conjecture. We hope to
study it in a future paper.


\section{Preliminary results}
\label{sec:notations}

\subsection{Formal group schemes}

We establish the notation and recall some of the results proved in
the papers \cite {Hurwitz} and \cite{GMP} that will be used in
subsequent sections.

In what follows, $\, V\,$ will denote a finite $\C((z))$--algebra
endowed with an action of the group $\Z/p \Z$, where $p$ is a
prime number, such that its fixed subset $V^{\Z/p\Z}$ is equal to
$\C((z))$. Let $\sigma$ denote a (fixed) generator of this
subgroup of $\aut_{\C((z))}V$.

Since $V$ is a finite $\C((z))$-algebra, there are canonical maps
given by the trace and the norm
    \begin{eqnarray*}
    \tr \colon V \, \longrightarrow & \,\C((z)) \\
    \nm \colon V \, \longrightarrow & \,\C((z))
    \end{eqnarray*}
which map an element $g$ in $V$ to the trace (respectively norm)
of the homothety of $V$ defined by $g$ (as a $\C((z))$-vector
space).

Throughout the paper, we will consider only the following two
cases.

\begin{enumerate}
 \item[(a)] \textsl{Ramified case:} $V = V_{\mathrm R}=\C((z_1))$, where the
$\C((z))$--algebra structure is given by mapping $z$ to $z_1^p$
and $\sigma$ is given by $\sigma(z_1)=\omega z_1$,  where $\omega$
is a primitive $p$-th root of $1$ in $\C$. In this case we set
$V^+_{\mathrm R}=\C[[z_1]]$ and $V^-_{\mathrm
R}=z_1^{-1}\C[z_1^{-1}]$.

\medskip

 \item[(b)] \textsl{Non-ramified case:}
$V = V_{\mathrm{NR}}=\C((z_1))\times \dots\times \C((z_p))$, where
the $\C((z))$--algebra structure is given by mapping $z$ to
$(z_1,\dots, z_p)$ and $\sigma$ is given by $\sigma(z_i)=z_{i+1}$
(for $i<p$) and $\sigma(z_p)=z_1$. In this case we set
$V^+_{\mathrm{NR}}=\C[[z_1]]\times \dots\times \C[[z_p]]$ and
$V^-_{\mathrm{NR}}=z_1^{-1}\C[z_1^{-1}]\times \dots\times
z_p^{-1}\C[z_p^{-1}]$.
\end{enumerate}

\begin{exam}
Let $C$ be a projective irreducible smooth curve with an order $p$
automorphism, denoted by $\sigma_C : C\to C$.

If $\pi : C\to \overline{C}:= C/<\!\sigma_C\!>$ denotes the
quotient map and $p\in \overline{C}$ is a smooth point, then the
$\w\o_{\overline{C},p}$-algebra $\w\o_{C,\pi^{-1}(p)}$ is
isomorphic to one of the above types.
\end{exam}

The formal base curve is $\w C:=\sf \C[[z]]$ and the formal
spectral cover is  $\w C_V:=\sf V^+$. Let $\Gamma_V$ be the formal
group scheme representing the functor
\begin{eqnarray*}
\left\{  \mbox{category of formal $\C$-schemes} \right\} \, &
\rightsquigarrow & \, \left\{
\text{category of groups} \right\} \\
S\, & \rightsquigarrow & \,  \left(V\hat{\otimes}_{\mathbb{C}}
H^0(S,\o_S)\right)^*_0
\end{eqnarray*}

\noindent where the subscript $0$ denotes the connected component
of the identity and the superscript ${}^*$ denotes the invertible
elements. Replacing $V$ by $V^+$ (respectively by $1+V^-$) we
define the subgroup $\bar \Gamma_V^+$ (respectively $\Gamma_V^-$)
and thus obtain the decomposition
    $$\Gamma_V=\Gamma_V^-\times \bar \Gamma_V^+ \ .$$
The formal Jacobian of the formal spectral cover is the formal
group scheme $\j(\w C_V):=\Gamma_V^-$. A straightforward
calculation shows that $\j(\w C_V)$ is the formal spectrum of the
ring
     {\small
    $$
    \o(\j(\w C_V))\,=\,
    \C\{\{t_1^{(1)},t_2^{(1)},\dots\}\}\w\otimes \dots \w\otimes
    \C\{\{t_1^{(r)},t_2^{(r)},\dots\}\}
    $$}
where $t_i^{(j)}$ are indeterminates, $\C\{\{t_1,t_2,\dots\}\}$
denotes the inverse limit $\limp \C[[t_1,\dots,t_n]]$, and $r=1$
for the ramified case (in which case the superscript $(1)$ is
dropped) and $r=p$ for the non-ramified one.

Replacing $V$, $V^+$ and $V^-$ by $\C((z))$, $\C[[z]]$ and
$z^{-1}\C[z^{-1}]$ respectively in the previous constructions, one
obtains formal schemes $\Gamma$, $\bar \Gamma^+$ and $\j(\w
C):=\Gamma^-$. It is straightforward that the canonical morphism
$\C((z))\hookrightarrow V$ gives rise to
$\Gamma\hookrightarrow\Gamma_V$ and that the trace and the norm
yield corresponding morphisms $\Gamma_V\to \Gamma$.

Recall that the Abel map $\phi_V:\w C_V\longrightarrow \j(\w C_V)$
is the morphism corresponding to the $\w C_V$-valued point of
$\Gamma_V$ associated to the $r$-tuple of series
    \begin{equation}\label{eq:formalAbelmorphism}
    \left( \left(1-\frac{\bar z_1}{z_1}\right)^{-1}, \dots,
    \left(1-\frac{\bar z_r}{z_r}\right)^{-1}\right)
    \end{equation}
where $\w C_V \simeq\sf \left( \C[[\bar z_1]]\times\dots \times
\C[[\bar z_r]]\right)$, with $r=1$ in the ramified case and $r=p$
in the non ramified one.

\begin{prop}\label{prop:albanese}
The Albanese variety of $\w C_V$ is the pair
$(\j(\w C_V),\phi_V)$.
\end{prop}

\begin{proof}
This statement is proved in~\cite{MP} for the case $V=\C((z))$.
The present case follows from that result and from the fact that
the Albanese variety of a disjoint union is the product of the
corresponding Albanese varieties.
\end{proof}

The natural morphism $\w C_V \stackrel{\pi}{\longrightarrow} \w C$
induces two group homomorphisms; namely, the pull-back
    $$ \pi^* \colon \j(\w C)\longrightarrow  \j(\w C_V)$$
and the Albanese
    $$\j(\w C_V)\longrightarrow \j(\w C) \ .$$

The Albanese map coincides with the restriction of the $\nm$ to
$\j(\w C_V)$; it will also be denoted by $\nm$.

Observe that the automorphism $\sigma$ of $V$ gives rise to
automorphisms of $\Gamma_V$ and $\j(\w C_V)$ (also denoted by
$\sigma$) such that $\Gamma_V^\sigma=\C((z))$ and $\j(\w
C_V)^\sigma=\j(\w C)$. In particular, it follows that $\tr$
(respectively $\nm$) maps an element $g\in \Gamma_V$ to
$\sum_{i=0}^{p-1}\sigma^i(g)$ (respectively
$\prod_{i=0}^{p-1}\sigma^i(g)$).

\subsection{Infinite Grassmannians}

Recall that the infinite Grassmannian $\grv$ of the pair $(V,V^+)$
is a $\C$-scheme, which is not of finite type, whose set of
rational points is
   {\small $$\left\{\begin{gathered}
   \text{subspaces $U\subset V$ such that }U\to V/V^+
   \\
   \text{has finite dimensional kernel and cokernel}
   \end{gathered}
   \right\} \ . $$}
This scheme is  equipped with the determinant bundle, $\Det_V$,
which is the determinant of the complex of $\o_{\grv}$-modules
$$
\L\,\longrightarrow\, V/V^+\,\hat\otimes_{\C}\, \o_{\grv} \, ,
$$
where $\L$ is the universal submodule of $\grv$ and the morphism
is the natural projection. The connected components of the
Grassmannian are indexed by the Euler--Poincar\'e characteristic
of the complex. The connected component of index $m$ will be
denoted by $\gr^m(V)$.

The group $\Gamma_V$ acts by homotheties on $V$, and this action
gives rise to a natural action on $\grv$
    $$\Gamma_V\times\grv\,\longrightarrow\,\grv \ .$$
Furthermore, this action preserves the determinant bundle.

These facts allow us to introduce $\tau$-functions and
Baker-Akhiezer functions of points of $\grv$. Let us recall the
definition and some properties of these functions
(\cite{Hurwitz},~\S3).

The determinant of the morphism $\L\to V/V^+\hat\otimes_{\C}
\o_{\grv}$ gives rise to a canonical global section
    $$
    \Omega_+\,\in\, H^0(\gr^0(V),\Detd_V) \ .
    $$
In order to extend this section to $\grv$ (in a non-trivial way),
we fix elements $v_m\in V^+$ for $m>0$ such that $\dim V^+/v_m
V^+=m$. Setting $v_{-m}:=v_m^{-1}$ for $m>0$, we define
$\Omega_+(U):=\Omega_+(v_m^{-1}U)$ for $U\in\gr^m(V)$.

Now, the $\tau$-function and BA functions will be introduced
following~\cite{Hurwitz}.  The $\tau$-function of $U$,
$\tau_U(t)$, is a function on $\j (\w C_V)$ and is introduced as a
suitable trivialization of the function $g\mapsto \Omega_+(g U)$
for $g\in \j(\w C_V)$,
    $$ \tau_U(g)
    \,=\, \frac{\Omega_+(gU)}{g \delta_U}
    $$
where $\delta_U$ is a non-zero element in the fibre of $\Detd_V$
over $U$.

Let $t$ be the set of variables $(t^{(1)},\dots, t^{(r)})$ (where
$t^{(j)}=(t_1^{(j)}, t_2^{(j)}, \dots)$) and $z_\centerdot$ denote
$(z_1,\dots,z_r)$. For $1\leq u, v \leq r$, let
$[z_v]:=(z_v,\frac{z_v^2}2,\frac{z_v^3}3,\dots)$,
$t+[z_v]:=(t^{(1)},\dots, t^{(v)}+[z_v],\dots, t^{(r)})$, $U_{uu}=
U$ and, if $u\neq v$,   $U_{uv}:=
(1,\dots,z_u,\dots,z_v^{-1},\dots,1)\cdot U$.

The $u$-th Baker-Akhiezer function of a point $U\in\grv$ is the
$V$-valued function
    {\small
    $$
    \psi_{u,U}({z_\centerdot},t)\,:=\, \left(\exp\Big(-\sum_{i\geq 1}
    \frac{t_i^{(1)}}{z_1^i}\Big)
    \frac{\tau_{U_{u1}}(t+[z_1])}{\tau_U(t)} , \dots, \exp
    \Big(-\sum_{i\geq 1} \frac{t_i^{(r)}}{z_r^i}\Big)
    \frac{\tau_{U_{ur}}(t+[z_r])}{\tau_U(t)}
        \right)
        \ .$$}

From the decomposition $V=\prod_{i=1}^r \C((z_i))$ we can write
    {\small
    $$\psi_{u,U}({z_\centerdot},t)\,=\,
    (\psi_{u,U}^{(1)}(z_1,t),\dots,\psi_{u,U}^{(r)}(z_r,t))
    \,\in\, \prod_{i=1}^r \left(
     \C((z_i))\w\otimes_\C \o(\j(\w C_V))\right) \ .$$}

Let $p_n(t)$ be the Schur polynomials and
$\wtilde\partial_{t^{(j)}} = (\partial_{t^{(j)}_1}, \frac{1}{2}
\partial_{t^{(j)}_2}, \frac{1}{3} \partial_{t^{(j)}_3}, \ldots )$.
Using the identity $\exp\left(\sum_{i\geq 1} z_j^i
\wtilde\partial_{t_i^{(j)}}\right)\tau_{U_{uj}}(t) =
\tau_{U_{uj}}(t+[z_j])$, we obtain
    {\small \begin{equation} \label{eq:BAjexpanssion}
    \begin{aligned}
    \psi_{u,U}^{(j)}(z_j,t)\,& =\, \exp\left(-\sum_{i\geq
    1}\frac{t_i^{(j)}}{z_j^i}\right) \frac{\tau_{U_{uj}}(t+[z_j])}{\tau_U(t)} \, =\\
    & =\, \left(\sum_{i\geq 0}\frac{p_i(-t^{(j)})}{z_j^i}\right)
    \frac{\left(\sum_{i\geq 0}
    p_i(\wtilde\partial_{t^{(j)}})z_j^i\right)\tau_{U_{uj}}(t)} {\tau_U(t)}\ .
    \end{aligned}
    \end{equation}}

The main property of these Baker-Akhiezer functions is that they
can be understood as generating functions for $U$ as a subspace of
$V$, as we recall next.

\begin{thm}[\cite{Hurwitz}] 
Let $U\in \gr^m(V)$. Then
    $$
    \psi_{u,U}({z_\centerdot},t)
    \,=\,
    v_m^{-1} \cdot (1,\dots,z_u,\dots,1) \cdot \sum_{i>0}\left(\psi_{u,U}^{(i,1)}(z_1), \dots,
    \psi_{u,U}^{(i,r)}(z_r)\right)p_{u i,U}(t)$$
where
    $$
    \big\{ (\psi_{u,U}^{(i,1)}(z_1), \dots,\psi_{u,U}^{(i,r)}(z_r)) \,\vert\,   i>0 \, , \, 1\leq u\leq r \big\}
    $$
is a basis of $U$ and $p_{u i,U}(t)$ are functions in $t$.
\end{thm}

Consider the following pairing
    $$
    \begin{aligned}
        V\times V &
    \,\longrightarrow\, \C \\
    (a,b) & \,\longmapsto \res_{z=0}\tr(a,b) \, dz \ .
    \end{aligned}
    $$
Since it is non-degenerate, there is an involution of $\grv$ which
maps any point $U$ to its orthogonal complement $U^\perp$. This
involution sends the connected component $\gr^m(V)$ to
$\gr^{1-m-p}(V)$ in the ramified case and to $\gr^{-m}(V)$ in the
non-ramified one.

Finally, the adjoint Baker-Akhiezer functions of $U$ are defined
by
    $$\psi_{u,U}^*({z_\centerdot},t)\,:=\, \psi_{u,U^\perp}({z_\centerdot},-t)\ ,$$
whose components are given by
    {\small \begin{equation}
    \label{eq:BAjexpanssion1}
    \psi_{v,U}^{*(j)}(z_j,s)\,=\, \left(\sum_{i\geq
    0}\frac{p_i(s^{(j)})}{z_j^i}\right) \frac{\left(\sum_{i\geq 0}
    p_i(-\wtilde\partial_{s^{(j)}})z_j^i\right)\tau_{U_{jv}}(s)} {\tau_U(s)}
     \ .
    \end{equation}}

\section{Moduli of curves with an order $p$ automorphism}
The aim of this section is to give explicit differential equations
that characterize the tau functions coming from curves with an
order $p$ automorphism among the tau functions coming from curves.

Recall from~\S5 of \cite{GMP} that $\M^\infty(p,{\mathrm R})$
(respectively  $\M^\infty(p,{\mathrm{NR}})$) is the subscheme of
$\gr(V_{\mathrm{R}})$ (respectively $\gr(V_{\mathrm{NR}})$)
parametrizing $(C,\sigma_C,x,t_x)$ where $C$ is a curve,
$\sigma_C$ is an order $p$ automorphism of $C$, $x$ is a smooth
point of $C$ fixed under $\sigma_C$ (respectively an orbit
consisting of $p$ pairwise distinct points) and $t_x$ is a formal
parameter at $x$.

More precisely, if $\M^\infty(1)$ (respectively $\M^\infty(p)$) is
the moduli space parameterizing $(C,x,t_x)$ where $C$ is a curve,
$x$ is a smooth point (respectively $p$ pairwise distinct smooth
points) and a formal parameter at $x$, then the Krichever map
induces an immersion
    $$
    \kr \colon \M^\infty(1) \,  \hookrightarrow\,
    \gr(V_{\mathrm R})
    $$
(respectively $\kr : \M^\infty(p) \hookrightarrow
\gr(V_{\mathrm{NR}})$) such that
    \begin{equation}
    \label{eq:MRequalMcutGrVR}
    \M^\infty(p,{\mathrm R}) \,  =\, \M^\infty(1) \cap \gr(V_{\mathrm
    R})^\sigma
    \end{equation}
(respectively $\M^\infty(p,{\mathrm{NR}}) = \M^\infty(p)
\cap\gr(V_{\mathrm{NR}})^\sigma $) where $\grv^\sigma$ denotes the
set of points in $\grv$ fixed under the action of $\sigma$.

Therefore, our task consists of writing down the hierarchies
corresponding to these subschemes.

\subsection{Hierarchies for invariant subspaces}

The automorphism of $\grv$ induced by $\sigma:V\to V$ will also be
denoted by $\sigma$; it preserves the determinant bundle. If we
denote by $\grv^\sigma$ the set of points in $\grv$ fixed under
the action of $\sigma$, then it is known that $\grv^\sigma$ is a
closed subscheme.

The action of $\sigma$ on $\grv$ can be easily described in terms
of the Baker-Akhiezer and the tau functions, as we show next.

Let us begin by studying the ramified case. We denote
    $$
    \lambda \, t \,=\, (\lambda \,  t_1, \lambda^2 t_2, \ldots ,
    \lambda^j t_j,\ldots) \text{ for any } \lambda\in \C
    $$
    and
    $$
    \sigma^*(t):=\omega^{-1} \, t \ .
    $$

Then
    $$
    \tau_{\sigma({U})}(t) \,=\, \tau_U(\omega^{-1} \, t)
    $$
(up to a constant) and
    \begin{equation}\label{eq:BAsigmaRam}
    \psi_{\sigma({U})}(z_1,t)
    \,=\, \psi_U(\omega^{-1} \, z_1, \omega^{-1}\, t)
    \,=\, \sigma^{-1}\left(\psi_U(\, z_1, \sigma^*(t))\right)
    \end{equation}
where the expression on the r.h.s. corresponds to the action of
$\sigma$ on $V$; that is, since $\psi_U(z_1,t)$ is $V$-valued,
$\sigma$ acts on $z_1$ and acts trivially on $t$.

It is also known that a point $U\in\grv$ lies in $\grv^\sigma$ if
and only if its BA functions satisfy the following identities
(\cite{Hurwitz}, Theorem~3.15).
    \begin{equation}
       \label{eq:characterizationgrsigma}
    \res_{z=0}\tr\left(\frac{1}{z_1}\psi_{\sigma(U)}(z_1,t)
    \psi_U^*(z_1,s)\right)\frac{dz}{z}\,=\,0
    \end{equation}

Then, one has the following

\begin{thm}[ramified case]\label{thm:equationsGrVR}
Let $V = V_{\mathrm{R}}$ and let $U$ be a closed point of $\grv$.

Then $U$ is a point of $\grv^\sigma$ if and only its
$\tau$-function $\tau_U$ satisfies  the following differential
equations
    {\small
    \begin{equation}
    \label{eq:characterizationgrsigmaram}
    \sum_{\substack{\beta_1 + \beta_2 - \alpha_1 -\alpha_2 = 1 \\
                    \alpha_1 -\beta_1 \equiv j  \pmod{p}}}
    \left( D_{\lambda_1  , \, \alpha_1}
    (-\tilde\partial_{x}) p_{\beta_1}(\tilde\partial_{x})
    D_{\lambda_2 , \, \alpha_2}(\tilde\partial_{s})
    p_{\beta_2}(-\tilde\partial_{s}) \right) \tau_U(x) \, \tau_U(s) = 0
    \end{equation}}
\noindent for all Young diagrams $\lambda_1,\lambda_2$ and all
integers $j\in \{0,\dots, p-1\}$, where the definition of
$D_{\lambda , \, \alpha}$ is
    $$
    D_{\lambda , \, \alpha}(\tilde\partial_{y})
    (f(y)) \,:= \, \sum_{\lambda -\mu = (\alpha)} \chi_{\mu}
    (\tilde\partial_{y}) (f(y))_{|_{y = 0}} \ .
    $$
\end{thm}

\begin{proof}
Let $U$ be a closed point of $\gr(V)$. The residue condition given
by (\ref{eq:characterizationgrsigma}) is equivalent to the
vanishing of the constant term of
    {\small
    $$ \sum_{j=1}^p \frac{1}{\omega^{j} z_1} \psi_{\sigma(U)}(\omega^j
     z_1 , t )
    \psi_U^*(\omega^j  z_1 , s)
    \,=\,
    \sum_{j=1}^p \frac{1}{\omega^{j} z_1} \psi_{U}(\omega^{j-1}  z_1
    , \omega^{-1}  t )
    \psi_U^*(\omega^j  z_1 , s)
    $$}

Using  (\ref{eq:BAjexpanssion}) and (\ref{eq:BAjexpanssion1}) this
coefficient turns out to be
    \begin{multline*}
    \sum_{j=1}^p \,
    \sum_{\substack{\beta_1 + \beta_2 - \alpha_1 -\alpha_2 = 1}}
    \omega^{j(\beta_2 + \beta_1 - \alpha_1 - \alpha_2
    -1) + \alpha_1-\beta_1}
    \cdot \\
    \cdot\, p_{\alpha_1}(- \omega^{-1} t)
    p_{\beta_1}(\tilde\partial_{\omega^{-1} t}) \tau_U(\omega^{-1}
    t) \, \cdot \,
    p_{\alpha_2}(s) p_{\beta_2}(-\tilde\partial_{s})\tau_U(s)
    \ .
    \end{multline*}

Since this coefficient must vanish for each $p$-th rooth
$\omega^j$ of $1$ (note that the case $\omega^0=1$ is precisely
the KP-hierarchy, see \cite{MP}), we obtain the equivalent
conditions
    $$
    \sum_{\substack{\beta_1 + \beta_2 - \alpha_1 -\alpha_2 = 1 \\
                    \alpha_1 -\beta_1 \equiv j \pmod{p} }}
    p_{\alpha_1}(- \omega^{-1} t)
    p_{\beta_1}(\tilde\partial_{\omega^{-1} t}) \tau_U(\omega^{-1} t)
    p_{\alpha_2}(s) p_{\beta_2}(-\tilde\partial_{s})\tau_U(s)
    \,=\,0
    $$
for all $j$ in $\{ 0, \ldots , p-1\}$ and all $t$ and $s$.

Equivalently (substituting $\omega^{-1} t$ by $x$), we obtain
    \begin{equation}
    \label{eq:fxs}
    \sum_{\substack{\beta_1 + \beta_2 - \alpha_1 -\alpha_2 = 1 \\
                    \alpha_1 -\beta_1 \equiv j \pmod{p} }}
    p_{\alpha_1} (-x) p_{\beta_1}(\tilde\partial_{x}) \tau_U(x)
    p_{\alpha_2}(s) p_{\beta_2}(-\tilde\partial_{s})\tau_U(s) \,=\, 0
    \end{equation}
\noindent for all values of $x$ and $s$, and for all $j$ in $\{ 0,
\ldots , p-1\}$.

Since the vanishing of a function $f(x,s)$ (such as the left hand
side of (\ref{eq:fxs})) for all values of $x$ and $s$ is
equivalent to the vanishing of $\chi_{\lambda_1}
(\tilde\partial_{x}) \chi_{\lambda_2} (\tilde\partial_{s})
f(x,s)_{|_{x=0, s=0}}$ for all Young diagrams $\lambda_1$ and
$\lambda_2$, a calculation shows that (\ref{eq:fxs}) is equivalent
to (\ref{eq:characterizationgrsigmaram}), thus proving the
theorem.
\end{proof}

\begin{rem}
We now compute the first equations in the previous statement and
relate them to the KP equations. Observe that if $E_j$ denotes the
l.h.s. of equation~(\ref{eq:characterizationgrsigmaram}) for
$j\in\{0,\dots,p-1\}$, then the KP hierarchy is $E_0+\dots+
E_{p-1}=0$. Let us make this explicit.

The first non trivial equation in the KP hierarchy, which
corresponds to the Young diagrams $\lambda_1 = (1,1,1)$ and
$\lambda_2 = 0$, is the celebrated KP equation
    \begin{multline*}
    \tau_U (0) \, \chi_{(2,2)}(\tilde\partial_{s}) \tau_U(s)_{|_{s =
    0}} - p_1(\tilde\partial_{x}) \tau_U(x)_{|_{x = 0}} \, \cdot
    \chi_{(2,1)}(\tilde\partial_{s}) \tau_U(s)_{|_{s = 0}}  \\+
    p_2(\tilde\partial_{x}) \tau_U(x)_{|_{x = 0}}
    \chi_{(1,1)}(\tilde\partial_{s}) \tau_U(s)_{|_{s = 0}}
    = 0 \ .
    \end{multline*}
On the other hand, equations~(\ref{eq:characterizationgrsigmaram})
for the same Young diagrams are
    {\small
    $$
    \left\{
    \begin{aligned}
    &
    p_1(\tilde\partial_{x}) \tau_U(x)_{|_{x =
    0}} \, \cdot \chi_{(2,1)}(\tilde\partial_{s}) \tau_U(s)_{|_{s =
    0}}  \,=\, 0
    \\
    &
    \tau_U (0) \, \chi_{(2,2)}(\tilde\partial_{s}) \tau_U(s)_{|_{s =
    0}} + p_2(\tilde\partial_{x}) \tau_U(x)_{|_{x = 0}}
    \chi_{(1,1)}(\tilde\partial_{s}) \tau_U(s)_{|_{s = 0}}
    \,=\, 0 \ ,
    \end{aligned}
    \right.
    $$}
for $p=2$ and
    {\small
    $$
    \left\{
    \begin{aligned}
    &
    p_1(\tilde\partial_{x}) \tau_U(x)_{|_{x =
    0}} \, \cdot \chi_{(2,1)}(\tilde\partial_{s}) \tau_U(s)_{|_{s =
    0}}  \,=\, 0
    \\
    &
    \tau_U (0) \, \chi_{(2,2)}(\tilde\partial_{s}) \tau_U(s)_{|_{s =
    0}} \,=\, 0
    \\
    &
    p_2(\tilde\partial_{x}) \tau_U(x)_{|_{x = 0}}
    \chi_{(1,1)}(\tilde\partial_{s}) \tau_U(s)_{|_{s = 0}}
    \,=\, 0 \ .
    \end{aligned}
    \right.
    $$}
for $p \neq 2$
\end{rem}

\begin{exam}
As an example we give some more equations
from~(\ref{eq:characterizationgrsigmaram}) for other pairs of
Young diagrams. Note that the corresponding equations in the KP
hierarchy are all trivial.

\begin{enumerate}
\item
Consider the Young diagrams  $\lambda_1 = 0$ and $\lambda_2 = (1)$
and the corresponding $p$ equations given by
(\ref{eq:characterizationgrsigmaram}).

If $p = 2$, both equations are trivial.

However, if $p \neq 2$, two of the equations
(\ref{eq:characterizationgrsigmaram}) (the cases $\alpha_1
-\beta_1 \equiv 0 \pmod p$ and $\alpha_1 -\beta_1 \equiv -2 \pmod
p$) are equivalent to the following
$$
\tau_U (0) \, \cdot \, p_2(\tilde\partial_{x}) \tau_U(x)_{|_{x =
0}} = 0 \ .
$$

Similarly, the consideration of $\lambda_1 = (1)$ and $\lambda_2 =
0$ in (\ref{eq:characterizationgrsigmaram}) yields trivial
equations for $p=2$, whereas  for $p \neq 2$ we obtain
$$
\tau_U (0) \, \cdot \, p_2(-\tilde\partial_{x}) \tau_U(x)_{|_{x =
0}}          = 0 \ .
$$

 \item
    More generally, for $n \geq 2$ and not divisible by $p$ we
    obtain
            \begin{align*}
        \tau_U (0) \, \cdot p_n(\tilde\partial_{t}) \tau_U(t)_{|_{t = 0}}
        & = 0 \\
        \tau_U (0) \, \cdot p_n(-\tilde\partial_{t}) \tau_U(t)_{|_{t = 0}} & = 0
        \end{align*}
     when considering (\ref{eq:characterizationgrsigmaram}) with
     $\lambda_i= (n-1)$ and $\lambda_j = 0$.
\item
    When $n \geq 2$ is not divisible by $p$, we obtain
        \begin{align*}
        p_1(\tilde\partial_{t}) \tau_U(t)_{|_{t = 0}} \, \cdot
        p_{n+1}(\tilde\partial_{s}) \tau_U(s)_{|_{s = 0}} & = 0 \\
        p_1(\tilde\partial_{t}) \tau_U(t)_{|_{t = 0}} \, \cdot
        p_{n+1}(-\tilde\partial_{s}) \tau_U(s)_{|_{s = 0}} & = 0
        \end{align*}
        by considering $\lambda_i= (n-1, 1)$ and $\lambda_j = 0$.
\item
    When $n \geq 1$ is not divisible by $p$, the equations become
        \begin{align*}
        p_2(\tilde\partial_{t}) \tau_U(t)_{|_{t = 0}} \, \cdot
        p_{n+2}(\tilde\partial_{s}) \tau_U(s)_{|_{s = 0}} & = 0 \\
        p_2(-\tilde\partial_{t}) \tau_U(t)_{|_{t = 0}} \, \cdot
        p_{n+2}(-\tilde\partial_{s}) \tau_U(s)_{|_{s = 0}} & = 0
        \end{align*}
    for $\lambda_i = (n+1, 2)$ and $\lambda_j = 0$.
\end{enumerate}
\end{exam}

\begin{rem}
Let us study the relation with the KdV hierarchy. Consider an
invariant point $U\in\grv^\sigma$ (in the ramified case). If we
impose the condition that $\C[z^{-1}]\cdot U =U$, then we obtain
the $p$-KdV hierarchy. To see this, recall that the
$\tau$-function of $U$, $\tau_U$, is (up to a constant) the
pullback of the global section $\Omega_+$ to $\Gamma_V$ by
    $$
    \Gamma_V\times \{U\}\,\longrightarrow\, \grv \ .
    $$

Since the condition means that $\Gamma^- \cdot U=U$, we obtain the
following diagram
    $$
    \xymatrix{
    {\mathcal P}
     \ar@^{(->}[r]
    \ar[rd]_\simeq &
    \Gamma_V^- \ar@{->>}[d] \ar[r] &  \grv  \\ &
    \Gamma_V^-/\Gamma^- \ar[ru] &
        }$$
where
    {\small
    $$
    {\mathcal P}\,
    :=\, \big\{g\in\j(\w C_V)\,\vert\, \nm(g)=1\big\}
    \,=\,
    \left\{ \exp\Big(\sum t_i z^{-i}_1\Big) \in \Gamma_V \,\vert \,
    t_i=0\text{ for } i= \dot{p}
    \right\}
    \ .
    $$}

Therefore, for the ramified case and for $p=2$, it makes sense to
write $\tau_U=\tau_U(t_1,t_3,\dots)$ as an element of
$\o({\mathcal P})=\C\{\!\{t_1,t_3,\dots\}\!\}$ so that the
$\tau$-function only depends on $t_i$ with $i$ odd. The resulting
hierarchy is the classical KdV hierarchy (see also \cite{SW},
Proposition~5.11).
\end{rem}

Now we focus in the non-ramified case. Denote
    $$
    \sigma^*(t)  \,:=\, (t^{(p)}, t^{(1)}, t^{(2)}, \ldots ,
    t^{(p-1)});
    $$
then the corresponding relations are as follows.
    {\small
    $$
    \begin{aligned}
    \tau_{\sigma({U})}(t)\, & =\,
    \tau_{\sigma({U})}((t^{(1)}, \ldots, t^{(p)}))
    \,=\\
    & =\, \tau_U((t^{(p)}, t^{(1)}, t^{(2)}, \ldots ,t^{(p-1)}))
    \,=\,
    \tau_U(\sigma^*(t))
    \end{aligned}
    $$}
(up to a constant) and
    {\small
    \begin{equation}\label{eq:sigmaofBAfunctionNonRam}
    \begin{aligned}
    \psi_{u,\sigma({U})}(z_{\centerdot} , t)
    \,& =\,
    (\psi_{u+1,U}^{(2)}(z_1, \sigma^*(t)),
    \psi_{u+1,U}^{(3)}(z_2,\sigma^*(t)),
    \ldots , \psi_{u+1,U}^{(1)}(z_p, \sigma^*(t)))\, =
    \\
    &=\,
    \sigma^{-1} \left(\psi_{u+1,U}(z_{\centerdot} ,
    \sigma^*(t))\right)
    \end{aligned}
    \end{equation}}
where the action of $\sigma$ on the r.h.s. is the action on
$V$-valued functions.

It is also known that a point $U\in\grv$ lies in $\grv^\sigma$ if
and only if its BA functions satisfy the following identities
(\cite{Hurwitz})
    \begin{equation}
    \label{eq:characterizationgrsigmaNR}
    \res_{z=0}\tr\left(
    \frac{\psi_{u,\sigma(U)}(z_\centerdot,t)}
    {(1,\dots,z_u,\dots,1)}
    \cdot
    \frac{\psi_{v,U}^*(z_\centerdot,s)}{(1, \dots , z_v ,\dots,1)}
    \right)
    {dz}\,=\,0
    \end{equation}
for all $u, v\in \{1,\dots,p\}$ where $(1,\dots,z_u,\dots, 1)$
denotes the element of $V$ with entries equal to $1$ except the
$u$-th, which is $z_u$.

In a similar manner to the ramified case  but this time using
(\ref{eq:characterizationgrsigmaNR}), we obtain the following
result.

\begin{thm}[non-ramified]\label{thn:sigmanr}
Let $V = V_{\mathrm{NR}}$ and let $U$ be a closed point of
$\gr(V)$.

Then $U$ is a point of $\grv^\sigma$ if and only the
$\tau$-functions of $U$, $\tau_{U_{uv}}$, satisfy the following
differential equations.
    {\small
    \begin{multline*}
    \sum_{j=1}^p \sum_{\substack{\beta_1+\beta_2-\alpha_1-\alpha_2    \\
    = \, \delta_{ju}+\delta_{jv}-1}} \left( D_{\lambda_j , \,
    \alpha_1}
    (-\tilde\partial_{x^{(j+1)}}) p_{\beta_1}
    (\tilde\partial_{x^{(j+1)}}) \text{\v{D}}_{\lambda_{\,
    \centerdot}}^{j+1} (\tilde\partial_{x})  \, \tau_{U_{u+1,j+1}}(x) \, \cdot \right.
    \\
    \left. \cdot \, D_{\mu_j , \, \alpha_2}(\tilde\partial_{s^{(j)}})
        p_{\beta_2}(-\tilde\partial_{s^{(j)}})
    \text{\v{D}}_{\mu_{\, \centerdot}}^{j} (\tilde\partial_{s}) \, \tau_{U_{jv}} (s)\right)   = 0
    \end{multline*}}

\noindent for all Young diagrams $\lambda_1, \mu_1 , \ldots ,
\lambda_p , \mu_p$, where  the differential operator
${\text{\v{D}}}_{\mu_{\, \centerdot}}^{\; k}$ is given by
$$
\text{\v{D}}_{\mu_{\, \centerdot}}^{\; k} (\tilde\partial_{s})
(f(s))= \prod_{\ell \neq k}  \chi_{\mu_{\ell}}
(\tilde\partial_{s^{\ell}}) (f(y))_{|_{s^{\ell} = 0}} \ .
$$
\end{thm}

\begin{proof}
Let $U$ be a closed point of $\gr(V)$. The residue condition given
by (\ref{eq:characterizationgrsigmaNR}) is equivalent to the
vanishing of
    {\small
    \begin{multline*}
    \res \left( \sum_{j=1}^p z^{-\delta_{ju} - \delta_{jv}} \,
    \psi_{u,\sigma(U)}^{(j)} (z,t) \psi_{v,U}^{*(j)} (z,s) \right) \,
    dz \,= \\
    = \, \res \left( \sum_{j=1}^p z^{-\delta_{ju} - \delta_{jv}} \,
    \psi_{u+1,U}^{(j+1)} (z,\sigma^*(t)) \psi_{v,U}^{*(j)} (z,s)
    \right) \, dz
    \end{multline*}}

Using  (\ref{eq:BAjexpanssion}) and (\ref{eq:BAjexpanssion1}) and
denoting $x = \sigma^*(t)$, the coefficient of $z^{-1}$ in the
latter sum turns out to be
    \begin{multline}
    \label{eq:fxsnr}
    \sum_{j=1}^p \,
    \sum_{\substack{\beta_1 + \beta_2 - \alpha_1 -\alpha_2\\
     = \, \delta_{ju} + \delta_{jv} -1}}
    \, p_{\alpha_1}(- x^{(j+1})p_{\beta_1}(\tilde\partial_{x^{(j+1)}}) \tau_{U_{u+1,j+1}}(x)\, \cdot \\
    \, \cdot \,  p_{\alpha_2}(s^{(j)}) p_{\beta_2}(-\tilde\partial_{s^{(j)}})\tau_{U_{jv}}(s)
    \ .
    \end{multline}

Since the vanishing of the function $f(x,s)$ one given by
(\ref{eq:fxsnr}) for all values of $x$ and $s$ is equivalent to
the vanishing of
    $$
    \displaystyle \prod_{1 \leq a, b \leq p}
    \chi_{\lambda_a} (-\tilde\partial_{x^{(a)}}) \chi_{\mu_b}
    (\tilde\partial_{s^{(b)}}) f(x,s)_{|_{x=0, s=0}}
    $$
for all Young diagrams $\lambda_1, \mu_1, \ldots, \mu_p$, a
calculation shows that the vanishing of (\ref{eq:fxsnr}) is
equivalent to (\ref{eq:characterizationgrsigmaNR}), thus proving
the theorem.
\end{proof}

\subsection{Equations of the moduli space}

Recalling the relation~\ref{eq:MRequalMcutGrVR} and the
Theorem~\ref{thm:equationsGrVR} we can now write down the
characterization of $\M^\infty(p,{\mathrm R})$ in terms of
differential equations for the $\tau$-functions.

\begin{thm}[ramified case]\label{thm:equationsMram}
Let $V = V_{\mathrm{R}}$ and let $U$ be a closed point of
$\gr^m(V)$.

Then $U$ is a point of $\M^\infty(p,{\mathrm R})$ if and only if
its $\tau$-function $\tau_U(t)$ satisfies the following set of
differential equations, for all Young diagrams $\lambda_1$,
$\lambda_2$, $\lambda_3$, $\lambda$, and all integers $j$ in $\{
0, 1, \ldots , p-1\}$.

\begin{enumerate}
\item The equations (\ref{eq:characterizationgrsigmaram}):
    {\small
    $$\sum_{\substack{\beta_1 + \beta_2 - \alpha_1 -\alpha_2 = 1 \\
                    \alpha_1 -\beta_1 \equiv j  \pmod{p}}}
    \left( D_{\lambda_1  , \, \alpha_1}
    (-\tilde\partial_{x}) p_{\beta_1}(\tilde\partial_{x})
    D_{\lambda_2 , \, \alpha_2}(\tilde\partial_{s})
    p_{\beta_2}(-\tilde\partial_{s}) \right) \tau_U(x) \, \tau_U(s) = 0
    $$}

\item
    {\small $$
    \sum_{\beta  - \alpha = m} \left( D_{\lambda  , \,
    \alpha}(\tilde\partial_{t}) p_{\beta}(-\tilde\partial_{t}) \right)
    \tau_U(t) = 0
    $$}

\item
    {\small \begin{multline*}
    \sum_{\beta_1 + \beta_2 +\beta_3 - \alpha_1 -\alpha_2
    -\alpha_3=2-m} \left( D_{\lambda_1  , \, \alpha_1}
    (-\tilde\partial_{t}) p_{\beta_1}(\tilde\partial_{t}) D_{\lambda_2
    , \,
    \alpha_2}(-\tilde\partial_{s}) p_{\beta_2}(\tilde\partial_{s}) \,\cdot \right.\\
    \left.\cdot \, D_{\lambda_3 , \, \alpha_3}(\tilde\partial_{u})
    p_{\beta_3}(-\tilde\partial_{3})\right) \tau_U(x) \, \tau_U(s) \,
    \tau_U(t) = 0 \ .
    \end{multline*}}
\end{enumerate}
\end{thm}

Similarly, we can write down differential equations for
$\M^\infty(p,\mathrm{NR})$ for the non ramified case. In this case
the non-trivial curves with automorphisms appear when $m \leq 0$;
if $-m = qp+f$ with $0 \leq f < p$, then $v_m = z_1^{q+1} \cdot
\ldots  \cdot z_f^{q+1}  \cdot z_{f+1}^{q} \cdot \ldots \cdot
z_p^{q}$, and $v_{-m} = v_m^{-1}$.

\begin{thm}[non ramified case]\label{thm:equationsMNram}
Let $V = V_{\mathrm{NR}}$ and let $U$ be a closed point of
$\gr^m(V)$.

Then $U$ is a point of $\M^\infty(p,\mathrm{NR})$ if and only if
the following set of differential equations is satisfied by its
$\tau$-functions $\tau_{U_{uv}}(t)$, for all Young diagrams
$\lambda_\centerdot=\{\lambda_1,\dots,\lambda_p\}$,
$\mu_\centerdot=\{\mu_1,\dots,\mu_p\}$,
$\nu_\centerdot=\{\nu_1,\dots,\nu_p\}$ and all integers $u,v,w$ in
$\{ 1, 2, \ldots , p\}$.

\begin{enumerate}
\item The equations of Theorem~\ref{thn:sigmanr}
    {\small
    \begin{multline*}
    \sum_{j=1}^p \sum_{\substack{\beta_1+\beta_2-\alpha_1-\alpha_2    \\
    = \, \delta_{ju}+\delta_{jv}-1}} \left(
    D_{\lambda_j , \,  \alpha_1}(-\tilde\partial_{x^{(j+1)}})
    p_{\beta_1} (\tilde\partial_{x^{(j+1)}})
    \text{\v{D}}_{\lambda_{\,\centerdot}}^{j+1} (\tilde\partial_{x})
    \, \tau_{U_{u+1,j+1}}(x) \, \cdot \right.
    \\
    \left. \cdot \,
    D_{\mu_j , \, \alpha_2}(\tilde\partial_{s^{(j)}})
    p_{\beta_2}(-\tilde\partial_{s^{(j)}})
    \text{\v{D}}_{\mu_{\, \centerdot}}^{j} (\tilde\partial_{s})
    \, \tau_{U_{jv}} (s)\right)   = 0
    \end{multline*}}

\item
    {\small \begin{multline*}
    \sum_{j=1}^f \sum_{\substack{\beta  - \alpha \\
    = \, \delta_{ju}-2-q}} \left(
    D_{\lambda_j , \, \alpha}(\tilde\partial_{t^{(j)}})
    p_{\beta}(-\tilde\partial_{t^{(j)}})
    \text{\v{D}}_{\lambda_\centerdot}^{j} (\tilde\partial_{t})
    \right) \tau_{U_{ju}} (t)\,+ \\
    +\, \sum_{j=f+1}^p \sum_{\substack{\beta  - \alpha \\
    = \, \delta_{ju}-1-q}} \left(
    D_{\lambda_j , \, \alpha}(\tilde\partial_{t^{(j)}})
    p_{\beta}(-\tilde\partial_{t^{(j)}})
    \text{\v{D}}_{\lambda_\centerdot}^{j} (\tilde\partial_{t})
    \right)
    \tau_{U_{ju}}(t) = 0
    \end{multline*}}

\item
    {\small \begin{multline*}
    \sum_{j=1}^f \sum_{\substack{\beta_1 + \beta_2 +\beta_3 - \alpha_1
    -\alpha_2 -\alpha_3 \\
    =\, q+\delta_{ju}+\delta_{jv}+\delta_{jw}}} \left(
    D_{\lambda_j  , \, \alpha_1} (-\tilde\partial_{t^{(j)}})
    p_{\beta_1}(\tilde\partial_{t^{(j)}})
    \text{\v{D}}_{\lambda_\centerdot}^{j} (-\tilde\partial_{t})
    D_{\mu_j , \, \alpha_2}(-\tilde\partial_{s^{(j)}})
    p_{\beta_2}(\tilde\partial_{s^{(j)}})
    \text{\v{D}}_{\mu_\centerdot}^{j} (-\tilde\partial_{s})
    \right.\\
    \left.
    D_{\nu_j , \, \alpha_3}(\tilde\partial_{x^{(j)}})
    p_{\beta_3}(-\tilde\partial_{x^{(j)}})
    \text{\v{D}}_{\nu_\centerdot}^{j} (-\tilde\partial_{x})
    \right) \tau_{U_{uj}}(t) \,
    \tau_{U_{vj}}(s) \, \tau_{U_{jw}}(x)
    \,+ \\
    \, +
    \sum_{j=f+1}^p \sum_{\substack{\beta_1 + \beta_2 +\beta_3 -
    \alpha_1 -\alpha_2 -\alpha_3 \,=
    \\
    =\, q+\delta_{ju}+\delta_{jv}+\delta_{jw}-1}} \left(
    D_{\lambda_j  , \, \alpha_1} (-\tilde\partial_{t^{(j)}})
    p_{\beta_1}(\tilde\partial_{t^{(j)}})
    \text{\v{D}}_{\lambda_\centerdot}^{j} (-\tilde\partial_{t})
    D_{\mu_j , \, \alpha_2}(-\tilde\partial_{s^{(j)}})
    p_{\beta_2}(\tilde\partial_{s^{(j)}})
    \text{\v{D}}_{\mu_\centerdot}^{j} (-\tilde\partial_{s})
    \right.\\
    \left.
    D_{\nu_j , \, \alpha_3}(\tilde\partial_{x^{(j)}})
    p_{\beta_3}(-\tilde\partial_{x^{(j)}})
    \text{\v{D}}_{\nu_\centerdot}^{j} (-\tilde\partial_{x})
    \right)
    \tau_{U_{uj}}(t) \,
    \tau_{U_{vj}}(s) \, \tau_{U_{jw}}(x)
    = 0 \ .
    \end{multline*}}
\end{enumerate}
\end{thm}

\section{$\tau$-functions and theta functions of Jacobians}

\subsection{Geometrical meaning of \lq\lq formal\rq\rq  objects}
\label{sec:geommeaning}

This section aims at motivating the relation between
$\tau$-functions and theta functions of Jacobians. In particular,
$\tau$-functions attached to a Riemann surface with marked points
will be defined following the works of Fay, Krichever, Shiota and
Adler--Shiota--van Moerbecke (\cite{Fay,Krichever,Shiota,ASvM}).

Throughout this section, $C$ will be an integral complete curve
over $\mathbb{C}$ of genus $g$. Let $J_{g-1}(C)$ denote the scheme
parametrizing invertible sheaves of degree $g-1$. For the sake of
clarity, we will assume $C$ to be smooth although most of the
results established here hold in greater generality.

We let $r=1$ in the ramified case and $r=p$ in the non-ramified
case.

Let us fix data $(C,\bar x,t_{\bar x})$, where $\bar
x=\{x_1,\dots, x_r\}$ are $r$ pairwise distinct points of $C$ and
$t_{\bar x}$ is a collection of formal parameters $\{ t_{x_1},
\dots, t_{x_r}\}$ giving corresponding isomorphisms $t_{x_j} :
\w\o_{C,x_j}\overset{\sim}\to \C[[z_j]]$.

\begin{prop}\label{prop:infinitesimalAlbanese}
For each invertible sheaf $L\in J_{g-1}(C)$, there is a canonical
morphism $\gamma : \j(\w C_V) \to  J_{g-1}(C)$ such that the
diagram
    $$
    \xymatrix{
    \w C_V \ar@{^(->}[r]^-\alpha \ar[d]_{\phi_V} &
    C \ar[d]^{\phi_L} \\
    \j(\w C_V) \ar[r]_\gamma & J_{g-1}(C)}
    $$
is commutative. Here $\phi_V$ is the Abel map and $\phi_L$ sends a
point $x'\in C$ to $L(r\cdot x'-\bar x)$.
\end{prop}

\begin{proof}
The data $(C,\bar x,t_{\bar x})$ gives rise to a canonical
morphism:
    $$
    \o_C\,\hookrightarrow\,\w\o_{C,\bar
    x} \,\overset\sim\longrightarrow\, V^+ \, .
    $$
Then we define $\alpha$ to be the induced morphism  between the
corresponding schemes.

Because of the Albanese property given in
Proposition~\ref{prop:albanese}, the composition $\w C_V\to
J_{g-1}(C)$ factors through the Abel map $\phi_V:\w C_V\to \j(\w
C_V)$, thus defining $\gamma$.
\end{proof}

Let  $J_{g-1}^\infty(C,\bar x)$ be the scheme parametrizing pairs
$(L,\phi)$ where $L\in J_{g-1}(C)$ and $\phi: \w L_{\bar
x}\overset{\sim}\to \w\o_{C,\bar x}$. It carries a canonical
action of $\bar \Gamma_V^+$ since this group acts by homotheties
on the trivialization of $L$. Summing up, there is an exact
sequence of group schemes
    $$
    0 \,\to\,  \bar \Gamma_V^+ \,\to\,
    J_{g-1}^\infty(C,\bar x) \,\to\, J_{g-1}(C)
    \,\to\, 0 \ .
    $$

This sequence has also a formal counterpart
    $$
    0 \,\to\,  \bar \Gamma_V^+ \,\to\,
    \Gamma_V \,\to\, \j(\w C_V)
    \,\to\, 0 \ .
    $$

Since this last sequence splits, $\Gamma_V\simeq \j(\w C_V)\times
\bar \Gamma_V^+$ and every element $g\in \Gamma_V$ can be written
as $(g_-,g_+)$ with $g_- \in \j(\w C_V)$, $g_+ \in \bar
\Gamma_V^+$ and $g=g_-\cdot g_+$.

\begin{prop}
Let $(L,\phi)$ be a point in $J_{g-1}^\infty(C,\bar x)$. Then the
Abel maps $\phi_V$ and $\phi_L$ have canonical lifts to $\Gamma_V$
and $J_{g-1}^\infty(C,\bar x)$ respectively, and the map $\gamma$
has a lift
    $$
    \xymatrix{ \Gamma_V \ar[r]^-{\w\gamma} \ar[d] &
    J_{g-1}^\infty(C,\bar x)
    \ar[d] \\
    \j(\w C_V) \ar[r]_{\gamma} &
    J_{g-1}(C)}
    $$
compatible with those of $\phi_V$ and $\phi_L$.
\end{prop}

\begin{proof}
Let $E(u,v)$ denote the prime form of $C$ as a meromorphic
function on $C\times C$. Then the meromorphic function
$\frac{E(x,x_i)}{E(x,x')}$ has a zero at $x=x_i$ and a pole at
$x=x'$. Moreover, if $z_i$ is a coordinate at $x_i$ such that
$z_i(x_i)=0$ and $\bar z_i:=z_i(x')$, then the expansion of this
function at $x_i$ is the following series.
    $$
    t_{x_i}\left(
    \frac{E(x,x_i)}{E(x,x')}\right)\,\in \,
    \left(1-\frac{\bar z_i}{z_i}\right)^{-1}\cdot (1+z_i
    \C[[z_i]]) \, .
    $$

Therefore, the $r$-tuple consisting of the expansions
    $$
    t_{\bar x}\left(\frac{E(x,x_1)}{E(x,x')},\dots,
    \frac{E(x,x_r)}{E(x,x')}\right)
    $$
corresponds to a morphism $\w C_V\to \Gamma_V$ which lifts the
Abel morphism $\phi_V$ defined by the
$r$-tuple~(\ref{eq:formalAbelmorphism}).

To define the lift of $\phi_L$, it is enough to observe that if
the line bundle $L$ carries a formal trivialization $\phi:\w
L_{\bar x}\overset\sim\to \w \o_{C,\bar x}$, then $L(r\cdot
x'-\bar x)$ is canonically endowed with the trivialization given
by
    $$\phi\cdot \prod_{i=1}^r \frac{E(x,x_i)}{E(x,x')} \ .$$

Finally, the lift of $\gamma$ is defined by
    $$
    \begin{aligned}
    \w\gamma\colon
    \Gamma_V
    & \longrightarrow \, J_{g-1}^\infty(C,\bar x) \\
    g \,  & \longmapsto  \left(L\otimes L_{g_-}, \phi\cdot g \right)
    \end{aligned}$$
where $L_{g_-}$ is given as follows: let $D_i$ be a small disk
around $x_i$ such that $z_i$ defines a coordinate in $D_i$, then
$L_{g_-}$ consists of gluing the trivial bundles on $C-\bar x$ and
on $\overset\circ D_1,\dots ,\overset\circ D_r$ by the transition
functions $(g_-)_1,\dots,(g_-)_r$ (see \cite{SW}, Remark~6.8 and
\cite{Shiota}, Lemma~4).
\end{proof}

Recall that the Krichever map associated to $(C,\bar x,t_{\bar
x})$ is the morphism
    $$
    \begin{aligned}
    J_{g-1}^\infty(C,\bar x)  & \to \, \grv
    \\
    (L,\phi) \, & \mapsto (t_{\bar x}\circ \phi)\big( H^0(C-\bar
    x,L)\big)
    \end{aligned}
    \ .
    $$

Then, we obtain the following

\begin{thm}\label{thm:liftequalaction}
Let $\kr$ denote the Krichever map associated to $(C,\bar
x,t_{\bar x})$. For $(L,\phi)\in J_{g-1}^\infty(C,\bar x)$, let $U
= \kr(L,\phi) \in \grv$.

Then the composition
    $$
    \Gamma_V\,\overset{\w\gamma}\longrightarrow\, J_{g-1}^\infty(C,\bar x)
    \,\stackrel{\kr}{\hookrightarrow} \, \grv
    $$
coincides with the morphism $\mu_U: \Gamma_V \cong \Gamma_V \times
\{U\}\to \grv$ mapping $g$ to $g\cdot  U$.
\end{thm}

\begin{proof}
We have to check that $(\kr \circ \, \w\gamma)(g)=g\cdot U$, but
this is a straightforward consequence of the definition of
$\w\gamma$.
\end{proof}

In particular, we have obtained morphisms
    $$
    \xymatrix{
    \Gamma_V \ar[r]^-{\w\gamma}  &
    J_{g-1}^\infty(C,\bar x) \ar@{^(->}[r]^-{\kr} \ar[d]^{\pi} & \grv
    \\
     &
    J_{g-1}(C) }
    $$
As a straightforward consequence of the determinantal construction
of the theta divisor in $J_{g-1}(C)$ and of the determinant bundle
on $\grv$, it follows that
    $$
    \kr^*\Det_V^* \,\simeq \, \pi^*\o_{J_{g-1}}(\Theta)
    $$

\begin{rem}
It is easy to check that there is a canonical isomorphism of
formal schemes
    $$
    \Gamma_V(U)/\bar \Gamma_V^+\,\simeq\, J_{g-1}(C)^{\w\;}_L \ ,
    $$
where $J_{g-1}(C)^{\w\;}_L$ denotes the formal completion of
$J_{g-1}(C)$ at the point $L$. If $A_U = \{v\in V : v\cdot
U\subseteq U\}$ is the stabilizer of $U$,  then the above
isomorphism yields, at the level of tangent spaces,
    $$
    T_1\left(\Gamma_V(U)/\bar \Gamma_V^+\right) \,=\,
    V/(A_U+V^+) \,\simeq\, H^1(C,\o_C)
    $$
since $A_U = t_{\bar{x}} (H^0(C-\bar x,\o_C))$. This is related to
Mulase's characterization of Jacobian varieties (\cite{Mulase}).
\end{rem}

\subsection{$\tau$-function attached to a Jacobian}
\label{sec:taujac}

The relations just described between the determinant bundle and
the theta line bundle induce an explicit relation between
$\tau$-functions and theta functions, which we now study.

Let  $J(C)$ denote the Jacobian variety of $C$. Choose a
symplectic basis $\{ \alpha_1 , \ldots , \alpha_g, \beta_1 ,
\ldots , \beta_g \}$ for $H_1(C, \mathbb{Z})$ and let $\{ \omega_1
, \ldots , \omega_g \}$ denote the corresponding canonical basis
of holomorphic $1$-forms on $C$; that is,
    $$
    \int_{\alpha_i}\omega_j \,=\, \delta_{ij}
    \quad \text{ and } \quad \int_{\beta_i}
    \omega_j \,=\, \Omega_{ij}
    $$
where $\Omega = \left( \Omega_{ij} \right)$ is the period matrix
of $C$. Then, as a complex torus, we have
    $$
    J(C)\, =\, \mathbb{C}^g/ \mathbb{Z}^g + \Omega \,\mathbb{Z}^g
    \ .
    $$

We also recall that the Riemann theta function associated to
$\Omega$ is the quasi-periodic function on $\C^g$ (the universal
cover of $J(C)$) given by
    $$
    \theta(z) \,=\,  \sum_{m \in \mathbb{Z}^g} \exp \left( 2\pi i m^{t} \, z
    + \pi i m^{t} \Omega m \right) \ .
    $$
By the Riemann-Kempf Theorem, there is an identification of
$J(C)$ and $J_{g-1}(C)$ such that the zero divisor of $\theta(z)$
corresponds to the divisor $\Theta\subset J_{g-1}(C)$. Therefore,
a point $\xi\in \C^g$ gives rise to a invertible sheaf $L_{\xi}\in
J_{g-1}(C)$.

Given $(C, \bar x=\{x_1 , \ldots , x_r\} ,t_{\bar x})$ there is a
canonical map
    $$
    \begin{aligned}
    H^0(C,\Omega_{C})\,&
    \longrightarrow\, \Omega_{\w C_V}\simeq\C[[z_1]]dz_1\times \dots
    \times \C[[z_r]]dz_r
    \\
    \omega_i \,& \mapsto \,
    \left(\dots, (a_{i1}^{(j)} + a_{i2}^{(j)} z_j + a_{i3}^{(j)} z_j^2
    + \ldots )d z_j, \dots\right)
    \end{aligned}
    $$
which maps each $\omega_i$ to its local expansions at the points
$x_1,\dots,x_r$. Since we have chosen bases on both sides, for
each $j \in \{ 1 , \ldots , r\}$ we have a $g \times \infty$
matrix over $\mathbb{C}$, $A^{(j)}$, associated to the map
$H^0(C,\Omega_{C})\to \C[[z_j]]dz_j$.  Moreover, $\rank A^{(j)}=g$
for each $j \in \{ 1 , \ldots , r\}$.

\begin{rem}
The transpose of the above map is related to the map induced by
Proposition~\ref{prop:infinitesimalAlbanese} at the level of
tangent spaces.
\end{rem}

We will now define the $\tau$-function of $(C, x_1 , \ldots , x_r
, z_1 , \ldots , z_r , \xi)$, where   $\{ x_1 , \ldots , x_r \}$
are $r$ different points in $C$, with respective local holomorphic
coordinates $z_1 , \ldots ,z_r$, and  $\xi\in\C^g$ is a point in
the universal cover of $J(C)$.

For each $j \in \{ 1 , \ldots , r\}$ and each natural number $n$,
we let $\eta_n^{(j)}$ denote the normalized meromorphic $1$-form
on $C$ with a unique pole of order $n + 1$ at $x_j$ of the form
$d( z_j^{-(n+1)}) + O(1)$ and such that
    $$
    \int^x \eta_n^{(j)} \,=\,  z_j^{-n} + O(z_j) \qquad \text{ at } x = x_j
    \ .
    $$

We also consider the complex numbers $q_{nm}^{(j)}$ defined by the
following identities.
 $$ \int^x \eta_n^{(j)}\, =\, z_j^{-n} - 2 \sum_{m=1}^{\infty}
 q_{nm}^{(j)} \frac{z_j^m}{m} \qquad \text{ at } x = x_j \ .
 $$

If we let $ t$ be the $r$-tuple $(t^{(1)} , \ldots , t^{(r)})$
where each $t^{(j)}$ is a family of variables $(t_1^{(j)},
t_2^{(j)} , \ldots )$, then we define the following quadratic form
    $$ Q(t) \,=\, \sum_{n,m \geq 1} q_{nm}^{(1)} t_n^{(1)} t_m^{(1)} + \ldots +
    \sum_{n,m \geq 1} q_{nm}^{(r)} t_n^{(r)} t_m^{(r)}
    $$
and the point $A(t)$ of $J(C)$ with values in $\C\{\{t^{(1)} ,
\ldots , t^{(r)}\}\}$ given by
    $$
    A(t) \,=\,  A^{(1)} t^{(1)}  + \ldots +
    A^{(r)} t^{(r)} \ .
    $$

Finally, the $\tau$-function of $(C, x_1 , \ldots , x_r , z_1 ,
\ldots , z_r , \xi)$ is defined as follows
    \begin{equation}
    \label{eq:jactau}
    \tau(\xi ,t) \,:=\,     \exp (Q(t)) \, \theta (A(t) + \xi)
    \end{equation}
for $t = (t^{(1)}, \ldots , t^{(r)})$.

It is worth pointing out that standard calculations shows that the
$\tau$-function $\tau(\xi ,t)$ coincides with the $\tau$-function
$\tau_U(t)$ of the point $U$ defined through the Krichever map
(\cite{SW,ASvM}).

Now, let us introduce the associated BA function for the ramified
case ($r=1$) by
    \begin{equation}
    \label{eq:actualpsi} \psi_{\xi}(z,t)\,:=\,    \exp(-\sum_{i\geq
    1}\frac{t_i}{z^i}) \cdot
        \frac{\tau(\xi,t+[z])}{\tau(\xi , t)}
    \end{equation}
and the adjoint BA function by
    \begin{equation}
    \label{eq:actualpsistar}  \psi_{\xi}^*(z,t)\,:=\, \exp(\sum_{i\geq
    1}\frac{t_i}{z^i}) \cdot
        \frac{\tau(\xi,t-[z])}{\tau(\xi , t)}
    \end{equation}
where $t=(t_1, t_2, \ldots )$ and $t+[z] =(t_1 +z, t_2 +z^2/2,
\ldots, t_n +z^n/n , \dots )$.

Then it follows from \cite{DJKM,Fay} that this BA function
coincides with the BA function of $U$ (the corresponding point
under the Krichever map) and that the following equality holds
    \begin{equation}
    \label{eq:formalKP}
    \res_{z=0} \psi_{\xi}(z,t)\psi_{\xi}^*(z,s)
    \frac{dz}{z^2} = 0
    \end{equation}
for all $t$ and $s$.

In the non-ramified case ($r=p$), BA functions will be introduced
as follows. Let $u,v$ be two integers in $\{1,\dots,p\}$ and let
$T_{uv}$ denote the homothety of $V$ given by the element $(1,
\ldots , z_u , \ldots , z_v^{-1} , \ldots , 1)$ if $u\neq v$, and
the identity if $u=v$. Then it is easy to check that $T_{uv}$
induces an automorphism of each connected component of $\grv$ and
that these automorphisms lift canonically to automorphisms of the
determinant bundle, which will also be denoted by $T_{uv}$. In
particular, we have automorphisms of $H^0(\grv,\Detd)$.

Recalling the bosonization isomorphism and the fact that
$\o(\Gamma^+_V)^*\simeq \o(\Gamma^-_V)$, we obtain isomorphisms
    $$
    T_{uv}^*\colon \o(\Gamma^-_V)
    \overset{\sim}\longrightarrow \o(\Gamma^-_V)
    $$
and we define
    $$
    \tau_{uv}(\xi,t)\,: =\, T_{uv}^*(\tau(\xi,t))
    $$
where $\tau(\xi,t)$ is given by~(\ref{eq:jactau}). Note that
$\tau_{U_{uv}}$ coincides with $T_{uv}^*(\tau_U)$ for any
$U\in\grv$.

We will also consider the corresponding formal $u-$BA functions as
follows
    \begin{equation}\label{eq:defnBAu}
    \psi_{u,\xi} (z_{\cdot},t) \,: =\,
    \left(
    \psi_{u,\xi}^{(1)} (z_1,t), \ldots ,
    \psi_{u,\xi}^{(r)} (z_r,t)\right)
    \end{equation}
where
    \begin{equation}\label{eq:defnBAuj}
    \psi_{u,\xi}^{(j)} (z_j,t) \,: =\,
    \exp\left(-\sum_{i\geq 1}
    \frac{t_i^{(j)}}{z_j^i}\right)
    \frac{\tau_{u j}(\xi,t+[z_j])}{\tau(\xi,t)} \ ,
    \end{equation}
and the formal adjoint $u-$BA functions by
    \begin{equation}\label{eq:defnBAadju}
    \psi_{u,\xi}^{*} (z_{\cdot},t) \,: =\,
    \left(
    \psi_{u,\xi}^{*(1)} (z_1,t) \ldots ,
    \psi_{u,\xi}^{*(r)} (z_r,t)\right)
    \end{equation}
where
    \begin{equation}\label{eq:defnBAadjuj}
    \psi_{u,\xi}^{*(j)} (z_j,t) \,: =\,
    \exp\left(\sum_{i\geq 1}
    \frac{t_i^{(j)}}{z_j^i}\right) \frac{\tau_{j u} (\xi,t
    - [z_j])}{\tau(\xi,t)} \ .
    \end{equation}

Then, from~\cite{Hurwitz}, we obtain that the BA functions satisfy
the $(1,\overset{p}\dots, 1)$-KP-hierarchy
    {\small $$
    \res_{z=0}\tr\left(
    \frac{\psi_{u,U}(z_\centerdot,t)}
    {(1,\dots,z_u,\dots,1)}
    \cdot
    \frac{\psi_{v,U}^*(z_\centerdot,s)}{(1, \dots , z_v ,\dots,1)}
    \right)
    {dz}\,=\,0
    $$}
for all $u, v\in \{1,\dots,p\}$. Equivalently, in terms of
$\tau$-functions we have
    \begin{multline}
    \label{eq:1p1KP}
    \sum_{j=1}^p \,
    \sum_{\substack{\beta_1 + \beta_2 - \alpha_1 -\alpha_2\\
     = \, \delta_{ju} + \delta_{jv} -1}}
    \, p_{\alpha_1}(- t^{(j})p_{\beta_1}(\tilde\partial_{t^{(j)}}) \tau_{U_{u,j}}(t) \\
    \, \cdot \,  p_{\alpha_2}(s^{(j)})
    p_{\beta_2}(-\tilde\partial_{s^{(j)}})\tau_{U_{jv}}(s) = 0
    \ .
    \end{multline}
Moreover, it follows from~\cite{Hurwitz} that these hierarchies
characterize when the BA-functions (respectively $\tau$-function)
are the BA-functions (respectively $\tau$-function) of a point of
$\grv$.

\section{Geometric characterization of Jacobian varieties with automorphisms in terms of Sato Grassmannian}

In this section we characterize the points of the Sato
Grassmannian that arise from geometric data over an algebraic
curve with automorphisms via the Krichever construction. We prove
that these points are those whose orbit under the action of
$\Gamma_V$ is finite dimensional (up to the action of $\bar
\Gamma_V^+$). A related result has been established in~\cite{GMP}.
This type of characterizations dates back to the approach of
Mulase (\cite{Mulase}).

Let us denote by $\overline\pic(C)$ the moduli space of rank $1$
torsion free sheaves on $C$.

We say that $(C,\bar x,L)$ is \textsl{maximal} for a curve $C$, a
divisor $\bar x$ composed by $r$ pairwise distinct points
$x_1,\dots,x_r$ in $C$, and $L\in \overline\pic(C)$ when one of
the following equivalent condition holds (see page 38
of~\cite{SW})
\begin{itemize}
\item let $(C',\bar x',L')$ be another such triple, and suppose
there exists $\psi: C'\to C$ a birational morphism such that
$\psi(\bar x')=\bar x$ and $\psi_* L'\iso L$; then $\psi$ is an
isomorphism.

\item The canonical map $\o_C\to \End(L)$ is an isomorphism.
\end{itemize}

\medskip

We begin with a detailed study of the ramified case  and then
generalize the results to the non-ramified case.

\subsection{Ramified case}

For the sake of notation, in this subsection the subscript
$\mathrm{R}$ in $V_\mathrm{R}$, $V_\mathrm{R}^+$,
$\Gamma_{V_\mathrm{R}}$, $\dots$ will be omitted.

\begin{defn}
\label{def:functorJram} Let
$\underline{\overline\pic}^\infty(p,{\mathrm R})$ be the
contravariant functor from the category of $\C$-schemes to the
category of sets defined by
    $$
    S\, \rightsquigarrow\,\{(C,\sigma_C,x,t_x,L,\phi_x)\}
    $$
where
\begin{enumerate}
    \item
    $p_C:C\to S$ is a proper and flat morphism whose fibres are
    geometrically integral curves.
    \item
    $\sigma_C:C\to C$ is an order $p$ automorphism  (over $S$).
    \item
    $x:S\to C$ is a smooth section of $p_C$, fixed under $\sigma_C$,  such that the Cartier divisor
    $x(s)$ is a smooth point of $C_s:=p_C^{-1}(s)$ for all closed points $s\in
    S$.
    \item
    $t_x$ is an equivariant formal parameter along $x(S)$; that is,
    an equivariant isomorphism of $\o_S$-modules
    $t_x\colon \w \o_{C,{x(S)}}\iso\w V_S^+$.
    \item
    $L\in \overline\pic(C)$ satisfies that $(C_s,x(s),L\vert_{C_s})$ is
    maximal for all closed point $s\in S$.
    \item $\phi_x$ is a formal trivialization of $L$
    along $x(S)$; that is, an isomorphism $\phi_x : \w
    L_{x(S)}\simeq \w\o_{C,x(S)}$.
    \item
    $(C,\sigma_C,x,t_x,L,\phi_x)$ and
    $(C',\sigma_{C'},x',t_{x'},L',\phi_{x'})$ are said to be
    equivalent when there is an isomorphism of $S$-schemes
    $C\iso C'$ compatible with all the data.
\end{enumerate}
\end{defn}

The Krichever morphism for the functor
$\underline{\overline\pic}^\infty(p,{\mathrm R})$ is the morphism
of functors
    $$
    \kr \colon \underline{\overline\pic}^\infty(p,{\mathrm R})
    \,\longrightarrow\,
    \underline{\gr(V)}
    $$
which sends the $S$-valued point $(C,\sigma_C,x,t_x,L,\phi_x)$ to
the following submodule of $\w V_S:=V\w\otimes\o_S$
    $$
    (t_x\circ \phi_x)\left(\limil{m}
    (p_{C})_* L(m\cdot x)\right)\,\subset\, \w V_S \ .
    $$

\begin{thm}
The functor $\underline{\overline\pic}^\infty(p,{\mathrm R})$ is
representable by a subscheme ${\overline\pic}^\infty(p,{\mathrm
R})$ of $\grv$.
\end{thm}

\begin{proof}
Consider the morphism from
$\underline{\overline\pic}^\infty(p,{\mathrm R})$ to $\grv\times
\grv$ which sends the $S$-valued point
$(C,\sigma_C,x,t_x,L,\phi_x)$ to the following pair of submodules
    {\small
    $$
    \left( \; t_x\left(\limil{m}
    (p_{C})_* \o_C(m\cdot x)\right)
    \, , \,
    (t_x\circ \phi_x)\left(\limil{m}
    (p_{C})_* L(m\cdot x)\right)\;\right)
    \, \in \,
    \grv\times \grv
    $$}
where $p_C: C\times S\to C$ is the projection.

From the inverse construction of the Krichever map
(\cite{Krichever, SW}) one has that this map is injective and that
the image is contained in the set $Z$ of those pairs $({\mathcal
A},{\mathcal L})$ in $\grv\times \grv$ such that
    $$
    \o_S\subset {\mathcal A}
    \quad ,\quad
    {\mathcal A}\cdot{\mathcal A}\subseteq {\mathcal A}
    \quad ,\quad
    {\mathcal A}\cdot {\mathcal L}\subseteq {\mathcal L}
     \quad ,\quad
    \sigma({\mathcal A}) = {\mathcal A} \ .
    $$

Let us examine the maximality condition. For $({\mathcal
A},{\mathcal L})$ satisfying the above conditions, let
$A_{\mathcal L}$ denote the stabilizer of $\mathcal L$
    $$
    A_{\mathcal L}\,:=\,
    \{ v\in \w V_S\text{ such that } v\cdot {\mathcal L}\subseteq {\mathcal
    L}\} \ ,
    $$
and let $(C',\sigma_{C'},x',t_{x'},L',\phi_{x'})$ be the geometric
data defined by the pair $(A_{\mathcal L},{\mathcal L})$. Then the
inclusion ${\mathcal A}\subseteq A_{\mathcal L}$ gives rise to an
equivariant morphism of $S$-schemes $\psi:C'\to C$ such that
$\psi(x')=x$ and $\psi_*L'\simeq L$. The maximality condition says
that $\psi_s$ is an isomorphism for every closed point $s\in S$.
That is, $A_{\mathcal L}$ is a finite ${\mathcal A}$-module such
that ${\mathcal A}_s = (A_{\mathcal L})_s$ for all $s$. Therefore
we have that ${\mathcal A}= A_{\mathcal L}$. Summing up, we are
interested on the subset $Z_0$ of $Z$ consisting of those pairs
$({\mathcal A},{\mathcal L})$ such that ${\mathcal A}= {\mathcal
A}_{\mathcal L}$.

From the proof of Theorem~6.5 of~\cite{MP} we know that the
condition ${\mathcal A}_{\mathcal L}\subseteq {\mathcal A}$, and
hence $Z_0$, is closed. The Krichever construction implies that
$Z_0$ represents $\underline{\overline\pic}^\infty(p,{\mathrm
R})$. Finally, $p_2\vert_{Z_0}:Z_0 \hookrightarrow \grv$ is a
closed immersion (where $p_2$ is the projection onto the second
factor), and the theorem is proved.
\end{proof}

Let us consider the following action of $\Gamma_V$ on $\grv^p:=
\grv\times\overset{p}\dots\times\grv$
    $$
    \begin{aligned}
    \mu^p\colon\Gamma_V\times \grv^p\,& \longrightarrow\, \grv^p
    \\
    (g,(U_1,\dots,U_p)) & \mapsto (g U_1,\dots, g U_p) \ .
    \end{aligned}
    $$

Then the morphism
   $$
    \begin{aligned}
    \grv\,& \hookrightarrow\,\grv^p
    \\
    U\, &\mapsto U_{\sigma}:=(U,\sigma(U),\dots,\sigma^{p-1}(U))
    \end{aligned}
    $$
\noindent is a $\Gamma_V$-equivariant closed immersion.

The orbit of $U_\sigma\in \grv^p$  under the action of $\Gamma_V$
is the schematic image of $\mu^p_{U_\sigma} := \mu^p
\vert_{\Gamma_V \times\{U_\sigma\}}$; it will be denoted by
$\Gamma_V({U_\sigma})$. Section 4 of~\cite{GMP} implies that
$\Gamma_V({U_\sigma})/\bar \Gamma_V^+$ is a formal scheme whose
tangent space is
    $$
    T_{U_\sigma}\left(\Gamma_V({U_\sigma})/\bar \Gamma_V^+\right)
    \,\simeq\,
    T_1\Gamma_V/(\ker d\mu^p_{U_\sigma}+ T_1\bar \Gamma_V^+) \ ,
    $$
where
    $$
    d\mu^p_{U_\sigma}\colon T_1\Gamma_V\,\longrightarrow\,
    T_{U_\sigma}\grv^p
    $$
is the map induced by $\mu^p_{U_\sigma}$ on the respective tangent
spaces.

\begin{thm}[ramified]\label{thm:algebraizableifffinitedimension}
Let $U$ be a closed point of $\grv$. Then the following conditions
are equivalent.
\begin{enumerate}
    \item $\dim_\C T_{U_\sigma}(\Gamma_V({U_\sigma})/\bar \Gamma_V^+) <\infty$
    where $V=V_{\mathrm R}$, and
    \item there exists $(C,\sigma_C,x,t_x,L,\phi_x)\in
    {\overline\pic}^\infty(p,{\mathrm R})$ such
    that its image by the Krichever morphism is $U$.
\end{enumerate}
\end{thm}

\begin{proof}
The result follows straightforwardly from similar arguments to
those of Theorems~4.13 and~4.15 of~\cite{GMP}.
\end{proof}

\begin{rem}\label{rem:finiteorbitRam}
If the conditions of the theorem hold, it follows that $\ker
d\mu^p_{U_\sigma} =\cap_{i=0}^{p-1} \sigma^i A_U$, with $A_U$ the
stabilizer of $U$. It is then straightforward to show that a
similar characterization exists in terms of two copies of $\grv$
instead of $p$ copies.
\end{rem}

\begin{thm}
Let $\underline\pic^\infty(p,{\mathrm R})$ denote the subfunctor
of $\underline{\overline\pic}^\infty(p,{\mathrm R})$ consisting of
those data $(C,\sigma_C,x,t_x,L,\phi_x)$ such that the fibres
$C_s$ are smooth curves for all closed points $s\in S$.

Then $\underline\pic^\infty(p,{\mathrm R})$ is representable by an
open subscheme $\pic^\infty(p,{\mathrm R})$ of
$\overline\pic^\infty(p,{\mathrm R})$.
\end{thm}

\begin{proof}
Let ${\mathcal C}\to {\overline\pic}^\infty(p,{\mathrm R})$ be the
universal curve. Then $\pic^\infty(p,{\mathrm R})$ is given by the
open subscheme of $\overline\pic^\infty(p,{\mathrm R})$ consisting
of the points $s$ such that ${\mathcal C}_s$ is a smooth curve.

\end{proof}

\begin{rem}
Observe that we have a forgetful map
    $$
    \xymatrix{
    \overline\pic^\infty(p,{\mathrm R})
    \ar@{->>}[r] &
    \M^\infty(p,{\mathrm R})}
    $$
which sends $(C,\sigma_C,x,t_x,L,\phi_x)$ to $(C,\sigma_C,x,t_x)$.
The fibre of $(C,\sigma_C,x,t_x)$ when $C$ is a smooth curve is
essentially the scheme studied in Section~\ref{sec:geommeaning}.
\end{rem}

\subsection{Non-ramified case}

Similarly to the ramified case we now consider the functor
$\underline{\overline\pic}^\infty(p,{\mathrm{NR}})$. To define it,
it suffices to replace conditions ($1$) and ($3$) in
Definition~\ref{def:functorJram}  respectively by
\begin{enumerate}
    \item[(1')] $p_C\colon C\to S$ is a proper and flat morphism whose
    fibres are geometrically reduced curves.
    \item[(3')] $x_i\colon S\to C$ are disjoint smooth sections
    of $p_S$ ($i=1,\dots , p$) such that $\sigma_C(x_i)=x_{i+1}$ for
    $i<p$ and $\sigma_C(x_p)=x_{1}$. Also, for every closed point
    $s\in S$ and each irreducible component of $C_s$, there
    is at least one $i$ such that $x_i(s)$ lies on that component.
\end{enumerate}
We also take $V$, $V^+$, $\Gamma_V$, $\dots$ to be
$V_{\mathrm{NR}}$, $V_{\mathrm{NR}}^+$,
$\Gamma_{V_{\mathrm{NR}}}$, $\dots$ in this case.

With arguments similar to those in the previous section, we can
prove the following results.

\begin{thm}
The functor $\underline{\overline\pic}^\infty(p,{\mathrm{NR}})$ is
representable by a subscheme
$\overline\pic^\infty(p,{\mathrm{NR}})$ of $\grv$.

The subfunctor $\underline\pic^\infty(p,{\mathrm NR})$ of
$\underline{\overline\pic}^\infty(p,{\mathrm{NR}})$ consisting of
those data $(C,\sigma_C,\bar x,t_x,L,\phi_x)$ such that the fibres
$C_s$ are smooth curves for all closed points $s\in S$ is
representable by an open subscheme $\pic^\infty(p,{\mathrm{NR}})$
of $\overline\pic^\infty(p,{\mathrm{NR}})$.
\end{thm}

\begin{thm}[non-ramified]\label{thm:finitedimNonRam}
Let $U$ be a closed point of $\grv$. Then the following conditions
are equivalent.
\begin{enumerate}
    \item $\dim_\C T_{U_\sigma}(\Gamma_V({U_\sigma})/\bar \Gamma_V^+) <\infty$  where $V=V_{\mathrm{NR}}$,
    \item there exists $(C,\sigma_C,\bar x,t_x,L,\phi_x)\in
    {\overline\pic}^\infty(p,{\mathrm{NR}})$ such
    that its image by the Krichever morphism is $U$.
\end{enumerate}
\end{thm}

\begin{rem}\label{rem:finiteorbitNRam}
Analogously to the ramified case, it follows that there is a
similar statement in terms of two copies of $\grv$ instead of $p$
copies.
\end{rem}

\section{Characterization of Jacobian theta functions of Riemann surfaces with non-trivial automorphisms}

In this section we give the conditions that a theta function of a
p.p.a.v. should satisfy in order to be the theta function of the
Jacobian of a smooth irreducible projective curve.

We begin with the proof of a generalization of the known theorems
of Mulase (\cite{Mulase}) and Shiota (\cite{Shiota}) in terms of
the Sato Grassmannian.

We will use the following notation.

Let $\Omega\in \C^g$ be a point in the Siegel upper half space
such that the principally polarized abelian variety $X_\Omega :=
\C^g/(\Z^g+\Omega\Z^g)$ is irreducible. Let $\theta(z) =
\theta(z,\Omega)$ denote the Riemann theta function of $X_\Omega$.

Let $r$ be a natural number, let $A^{(j)} =
(a_1^{(j)},a_2^{(j)},\cdots)\in (\C^g)^\infty$ be a $g \times
\infty$-matrix of rank $g$ for each $j\in\{1 , \ldots , r \}$, let
$Q^{(j)}(t^{(j)}) = \sum_{i,k=1}^\infty q_{ik}^{(j)} t_i^{(j)}
t_k^{(j)}$, with $q_{ik}^{(j)}  \in \C$, be a quadratic form for
each $j$ in $\{1 , \ldots , r \}$, and let $\xi$ be in $\C^g$.

\begin{defn}\label{defn:tauBAofppav}
The $\tau$-function and BA-functions associated to the data
$(X_{\Omega},\Omega,\{A^{(1)},\dots,A^{(r)}\},
\{Q^{(1)},\dots,Q^{(r)}\}, \xi)$ are given formally by
formulae~(\ref{eq:jactau})
and~(\ref{eq:defnBAu})-(\ref{eq:defnBAadjuj}).

The $\tau$-function will be denoted by $\tau(\xi,t)$ while the
BA-function (respectively adjoint BA-function) will be denoted by
$\psi_{u,\xi}(z_\centerdot,t)$ (respectively
$\psi^*_{u,\xi}(z_\centerdot,t)$).
\end{defn}

\begin{thm}\label{thm:S-Mgen}
Let $X_\Omega$ be an irreducible p.p.a.v. of dimension $g$.

Then the following conditions are equivalent.

\begin{enumerate}
    \item There exists a triple $(C, \bar x, t_{\bar x})$, where
    $C$ is a  projective irreducible smooth curve of genus $g$, $\bar{x} =(x_1 , \ldots , x_r)$
    is an $r$-tuple of distinct points in $C$ and $t_{\bar x} = (t_{x_1} , \ldots , t_{x_r})$
    is an $r$-tuple of local parameters at the corresponding
    $x_j$, such that $X_\Omega$ is isomorphic as a p.p.a.v. to the Jacobian of $C$.

    \item For each $j \in  \{ 1 , \ldots , r \}$  there exist a $g \times
    \infty$-matrix of rank $g$, $A^{(j)} =
    (a_1^{(j)},a_2^{(j)},\cdots)$,
    with $a_i^{(j)} \in \C^g$, and a quadratic form
    $Q^{(j)}(t^{(j)}) = \sum_{i,k=1}^\infty q_{ik}^{(j)} t_i^{(j)}
    t_k^{(j)}$, with $q_{ik}^{(j)}  \in \C$, such that for every $\xi\in \C^g$,
    the corresponding $\tau$-function $\tau(\xi,t)$
    is a $\tau$-function of the $(1, \stackrel{r}{\ldots} ,
    1)-$KP-hierarchy (\ref{eq:1p1KP}).
\end{enumerate}

    Moreover, if one of the conditions is fulfilled, the
    matrices $A^{(j)}$ and the quadratic forms $Q^{(j)}$ are
    equal to the data associated to the triple $(C, \bar x, t_{\bar x})$ in Section
    4.B.
\end{thm}

\begin{proof}
$1. \Rightarrow 2.$ It follows from Section 4.

$2. \Rightarrow 1.$ We denote
    $$
    A(t) = \sum_{j=1}^r A^{(j)} t^{(j)} \quad \text{ and }\quad Q(t) =
    \sum_{j=1}^r \sum_{i,k \geq 1} q_{ik}^{(j)} t_{i}^{(j)}
    t_{k}^{(j)} \ .
    $$

Since $\tau(\xi,t)$ is a $\tau$-function of the $(1,
\stackrel{r}{\ldots} , 1)-$KP-hierarchy for every $\xi\in \C^g$,
it follows that $\tau(\xi,t)$ defines a point $U_\xi\in \grv$
(with $V=\C((z))\times\stackrel{r}\dots\times\C((z))$) such that
$\tau(\xi,t) = \tau_{U_{\xi}}(t)$ (up to a constant). From Theorem
3.12 in \cite{Hurwitz}, we have that
    {\small $$
    U_{\xi} = \left\langle \left(
    p_i(\widetilde{\partial_t}){\psi_{u,\xi}^{(1)}
    }(z_1,t)_{|_{t=0}} , \ldots ,
    p_i(\widetilde{\partial_t}){\psi_{u,\xi}^{(r)}}
    (z_r,t)_{|_{t=0}}\right) , i \geq 0, 1 \leq u \leq r
    \right\rangle .
    $$}
Therefore, we have obtained a morphism
    $$
    \begin{aligned}
    \varphi\colon \C^g
        \, & \longrightarrow \, \grv  \\
        \xi \,  & \longmapsto  U_\xi \ .
    \end{aligned}
    $$

We claim that this morphism induces an injection
    \begin{equation}\label{eq:XintoGrV}
    X \hookrightarrow \grv \ .
    \end{equation}
Indeed, given $\xi_1$ and $\xi_2$ in $\C^g$, the condition
$U_{\xi_1} = U_{\xi_2}$ is equivalent to $\tau(\xi_1,t) =
\tau(\xi_2,t)$ for all $t$ (up to a constant), which is in turn
equivalent to $\theta(A(t) + \xi_1) = \theta(A(t) + \xi_2)$ for
all $t$ (up to a constant), and therefore equivalent to $\xi_1
-\xi_2 \in \mathbb{Z}^g + \Omega \mathbb{Z}^g$, since $\Theta$ is
a principal polarization on $X$.

Now the function $A$ can be interpreted as a surjective linear map
    $$
    \C^{\infty} \times \stackrel{r}{\ldots} \times \C^{\infty}
    \stackrel{A}{\longrightarrow} \C^g
    $$
and, with the identifications $\C^{\infty} \times
\stackrel{r}{\ldots} \times \C^{\infty} = T_1 \Gamma^-_V$ and
$T_{\xi} X = \C^g$, $A$ corresponds to a surjective morphism of
formal group schemes
    $$
    \Gamma^-_V \stackrel{A_{\xi}}{\longrightarrow} \hat{X}_{\xi} \, .
    $$

We claim now that the surjective morphism
    $$
    \begin{aligned}
    \mu_\xi\colon \Gamma_V^- & \longrightarrow \, \Gamma(U_\xi)/{\bar
    \Gamma_V^+}
    \\
    g \, &\longmapsto  g\cdot U_\xi
    \end{aligned}
    $$
factorizes by $A_{\xi}:\Gamma_V^-\to \hat{X}_\xi$. Observe that if
$s =(s^{(1)}, \ldots , s^{(r)}) \in \ker A$ then
    {\small $$
    \tau_{U_{\xi}}(t+s)=\tau (\xi,t+s) = q_{\xi} (t,s) \,
    \exp(Q(t))\, \theta(A(t)+\xi)= q_{\xi} (t,s) \,
    \tau(\xi,t)
    $$}
where $q_{\xi} (t,s)$ is an exponential of a linear function in
$t$. Generalizing Lemma~3.8 of~\cite{SW} to the case of
$\Gamma_V$, there exists $g \in \bar \Gamma^+_V$ (which depends on
$s$) such that
    $$
    \tau_{U_{\xi}}(t+s)=q_{\xi} (t,s) \, \tau(\xi,t) =
    \tau_{g\cdot U_{\xi}}(t) \ .
    $$

Hence there is a factorization
    $$
    \xymatrix{
    \Gamma^-_V \ar@{->>}[rr]^{\mu_{\xi}} \ar[rd]_{A_{\xi}} & &
    \Gamma(U_{\xi})/\bar\Gamma^+_V
    \\
    &\hat{X}_{\xi} \ar@{-->>}[ru]}
    $$

In particular, it follows that $\dim T_{U_{\xi}}
\Gamma(U_{\xi})/\bar\Gamma^+_V$ is finite, and applying the
results of~\cite{Mulase} one has that there exists data
$(C_\xi,\bar{x}_\xi, t_{\xi}, L_\xi , \phi_{\xi})$ associated to
$U_{\xi}$ under the Krichever map.

Let us check that the piece of data $(C_\xi,\bar{x}_\xi, t_{\xi})$
does not depend on $\xi$. Indeed, for $\xi,\xi'\in X$ let $s\in
\C^{\infty} \times \stackrel{r}{\ldots} \times \C^{\infty}$ be
such that $A(s)=\xi'-\xi$. Then
    $$
    \begin{aligned}
    \tau_\xi(t+s)&  \,=\, \exp(Q(t+s))\theta(A(t+s)+\xi)
    \,=
    \\
    & \,=\,
    q_{\xi'}(t)\exp(Q(t))\theta(A(t)+\xi')\,=\,
    q_{\xi'}(t)\tau_{\xi'}(t)
    \end{aligned}
    $$
The generalization of Lemma~3.8 of~\cite{SW} implies that there
exists $g_{\xi'}(z)\in V$ such that $U_{\xi'}=g_{\xi'}(z) U_\xi$.
From this fact it follows that $U_{\xi'}$ and $U_\xi$ have the
same stabilizer and that, therefore, $(C_\xi,\bar{x}_\xi,
t_{\xi})$ does not depend on $\xi$. It will be denoted by $(C,\bar
x,t_{\bar x})$.

The latter fact does have further consequences. It implies that
the map~(\ref{eq:XintoGrV}) takes values into
$\overline{\pic}^\infty(C,\bar x)$ (the subscheme of $\grv$
parameterizing torsion free sheaves of rank $1$ on $C$ with a
formal trivialization along $\bar x$). Furthermore, it says that
the composite map
    $$
    \xymatrix@R=5pt{
    X \ar@{^{(}->}[r] &
    \overline{\pic}^\infty(C,\bar x) \ar@{->>}[r] &
    \overline{\pic}(C)
    \\
    \xi' \ar@{|->}[r] & U_{\xi'} \ar@{|->}[r] & L_{\xi'}
    }$$
takes values in $\pic^0(C)\cdot L_{\xi}$, the orbit of $L_{\xi}\in
\overline{\pic}(C)$ under the action of $\pic^0(C)$. Using the
surjectivity of $A$ we can show that the induced map
    $$
    X \,\longrightarrow\, \pic^0(C) \cdot L_{\xi}
    $$
is surjective. Since $(C,\bar x,t_{\bar x}, L_{\xi},\phi_\xi)$ is
maximal (see~\S5), the action of $\pic^0(C)$ on
$\overline{\pic}(C)$ is free. So $\pic^0(C)$ is a quotient of an
abelian variety and, therefore, $C$ is a smooth complete curve of
genus at most $g$.

To finish the implication $2. \Rightarrow 1.$, one has only to
show that $X \to\pic^0(C)$  is an isomorphism of p.p.a.v.'s: Given
$(X,\xi)$ and $(J(C),\xi)$ we consider the tau-functions $\tau_X =
\tau (\xi,t)$ associated to $X$ as in (\ref{defn:tauBAofppav}) and
$\tau_J = \tau (\xi,t)$ associated to $J(C)$ as in
(\ref{eq:jactau}). By the construction of the data
$(C_\xi,\bar{x}_\xi, t_{\xi}, L_\xi , \phi_{\xi})$, it follows
that $\tau_X = \tau_J$ (up to a constant) and hence
    $$
    \theta_X (A (t) + \xi) =  \exp(q(t)) \, \theta_J(A_J(t) +
    \xi) \ \text{ (up to a constant)} ,
    $$
where $q(t)$ is a quadratic function.

Therefore
    $$
    \Theta_X = \Theta_{J} \ \text{ (up to translation)}.
    $$
In particular, the curve $C$ is irreducible and of genus $g$.
\end{proof}

\begin{rem}
Considering $r=1$ in the previous theorem we obtain the
characterization of Jacobian varieties given by Shiota
(\cite{Shiota}, Theorem 6).
\end{rem}

We will now apply this result to give a sufficient and necessary
condition for a theta function of a p.p.a.v. to be the theta
function of a curve with an automorphism of prime order $p$ with a
fixed point.

\begin{thm}[ramified case]\label{thm:autram}
Let $X_\Omega$ be an irreducible p.p.a.v. of dimension $g$.

Then the following conditions are equivalent.
\begin{enumerate}
\item There exists a quadruple $(C, \sigma_C, x, t_{x})$, where
$C$ is a  projective irreducible smooth curve of genus $g$,
$\sigma_C$ is an automorphism of order $p$ of $C$, $x$ is a fixed
point of $\sigma_C$ in $C$, and $t_{x}$ is a local parameter at
$x$, such that $X_\Omega$ is isomorphic as a p.p.a.v. to the
Jacobian of $C$.

\item There exist  a $g \times \infty$-matrix $A$ of rank $g$ and
a symmetric quadratic form $Q(t)$  such that for each $\xi_0$ in
$\C^{g}$ there exists $\xi_{1}$ in $\C^{g}$ so that the
corresponding BA-functions satisfy
    \begin{enumerate}
    \item[a)] the KP hierarchy
        $$
        \res_{z=0} \psi_{\xi_0}(z,t)
        \psi_{\xi_0}^* (z,s) \frac{dz}{z^2}\,=\,0 \ , \text{and}
        $$
    \item[b)] the identity
        $$
        \res_{z=0}
        \psi_{\xi_0}(\omega^{-1} z,\omega^{-1} t)
        \psi_{\xi_1}^*(z,s) \frac{dz}{z^2}\,=\,0
        $$
    \end{enumerate}
for all $t$ and $s$, where $\omega$ is a primitive $p$-th root of
1.
\end{enumerate}
\end{thm}

\begin{proof}
$1. \Rightarrow 2.$

From the results given in Section \ref{sec:taujac} with $r=1$, we
know that given  $(C, \sigma_C, x, t_{x})$ as in condition 1.,
there exist $A$ and $Q$ as in condition 2) such that its
associated BA-functions $\psi_{\xi_0}(z,t)$ and
$\psi_{\xi_0}^*(z,s)$, defined by (\ref{eq:actualpsi}) and
(\ref{eq:actualpsistar}) respectively, satisfy (\ref{eq:formalKP})
for each $\xi_0$ in $\C^g$. That is, there exists a point
$U_{\xi_0}\in\gr(V_{\mathrm R})$ such that
$\psi_{\xi_0}=\psi_{U_{\xi_0}}$ and 2.a) follows.

Furthermore, the induced embedding $J(C)\hookrightarrow
\gr(V_{\mathrm R})$ (given by~(\ref{eq:XintoGrV})) is compatible
with the actions of $\sigma_C$ in $J(C)$ and of $\sigma$ on
$\gr(V_{\mathrm R})$; that is, $\sigma(U_\xi)=U_{\sigma_C^*\xi}$.
Having this in mind and taking into account the relation between
the BA-function of $U_{\xi_0}$ and $U_{\xi_1}=\sigma(U_{\xi_0})$
given by~(\ref{eq:BAsigmaRam}), 2.b) follows.

$2. \Rightarrow 1.$

We know by Theorem \ref{thm:S-Mgen} (with $r=1$) that the first
identity implies that there exists a triple $(C,x,t_x)$ such that
$X_{\Omega} \cong J(C)$ as p.p.a.v.'s.

The second identity implies that $U_{\xi_1} = \sigma(U_{\xi_0})$.
In particular, we have that $A_{\xi_1} = \sigma(A_{\xi_0})$ where
$A_\xi$ denotes the stabilizer of $U_\xi$. Then the orbit
$\Gamma_{V_{\mathrm R}}(U_{\xi_0}, \sigma(U_{\xi_0})/\bar
\Gamma_{V_{\mathrm R}}^+$ is finite dimensional and, by
Theorem~\ref{thm:algebraizableifffinitedimension} (see
Remark~\ref{rem:finiteorbitRam}), we obtain that $\sigma$ induces
an automorphism of $A_{\xi_0}$. Since $A_{\xi_0}=t_{\bar
x}(H^0(C-\bar x,\o_C))$, because of the Krichever construction,
the result follows.
\end{proof}

\begin{rem}
Under the hypotheses of the above Theorem, suppose there are
$\xi_0$ and $\xi_1$ in $\C^g$ such that $\xi_0 -\xi_1 \in
\mathbb{Z}^g + \Omega\mathbb{Z}^g $ and the equations of the
theorem are satisfied. Then there exists a line bundle $L$ on $C$
such that $\sigma^*(L) \simeq L$.
\end{rem}

\begin{thm}[non-ramified]\label{thm:autNram}
Let $X_\Omega$ be an irreducible p.p.a.v. of dimension $g$.

Then the following conditions are equivalent.
\begin{enumerate}
    \item
    There exists a quadruple $(C, \sigma_C, \bar{x}, t_{\bar{x}})$, where
    $C$ is a  projective irreducible smooth curve of genus $g$, $\sigma_C$ is
    an automorphism of order $p$ of $C$, $\bar{x} = \{ x_1 , \ldots , x_p\}$
    is an orbit of $\sigma_C$ consisting of $p$ different points in $C$, and
    $t_{\bar{x}} = \{ t_{x_1} , \ldots , t_{x_p} \}$ is a collection
    of local parameters $t_{x_j}$ at each respective $x_j$, such that $X_\Omega$ is
    isomorphic as a p.p.a.v. to the Jacobian of $C$.
    \item
    There exist $p$ matrices $A^{(1)},\dots, A^{(p)}$, where $A^{(j)}$
    is a $g \times \infty$-matrix of rank $g$, and $p$ symmetric quadratic
    forms $Q^{(1)},\dots, Q^{(p)}$, such that for each $\xi_0\in \C^g$ there exists
    $\xi_1\in \C^g$ so that their BA-functions satisfy
        \begin{enumerate}
        \item[a)] the $(1,\stackrel{r}{\ldots}, 1)$-KP
        hierarchy
            $$
            \res \left( \sum_{j=1}^p z^{-\delta_{ju} - \delta_{jv}} \,
            \psi_{u,\xi_0}^{(j)} (z,t) \psi_{v,\xi_0}^{*(j)} (z,s) \right) \, dz
            \,=\,0 \ , \text{and}
            $$
        \item[b)] the identity
            $$
            \res \left( \sum_{j=1}^p z^{-\delta_{ju} - \delta_{jv}} \,
            \psi_{v+1,\xi_0}^{(j+1)} (z,\sigma^*(t))
            \psi_{u,\xi_1}^{*(j)} (z,s) \right) \, dz
            \,=\,0
            $$
            where $\sigma^*(t):=
            (t^{(p)}, t^{(1)}, t^{(2)}, \ldots , t^{(p-1)})$.
        \end{enumerate}
\end{enumerate}
\end{thm}

\begin{proof}
$1. \Rightarrow 2.$

Given $(C, \sigma_C, \bar{x}, t_{\bar{x}})$ satisfying condition
1., we construct $A(t)$ and $Q(t)$ as in Section \ref{sec:taujac}.
By Theorem \ref{thm:S-Mgen}, for each $\xi \in \C^g$ the
corresponding BA-functions satisfy 2.a).

Similarly to the ramified case, the actions of $\sigma_C$ in
$J(C)$ and of $\sigma$ on $\gr(V_{\mathrm{NR}})$ are compatible;
that is, $\sigma(U_\xi)=U_{\sigma_C^*\xi}$. Taking into account
the relation~(\ref{eq:sigmaofBAfunctionNonRam}) between
$\psi_{u,\xi_0}$ and $\psi_{v,\xi_1}$ for
$L_{\xi_1}=\sigma_C^*L_{\xi_0}$, the second part of 2. follows.

$2. \Rightarrow 1.$

From 2.a) and Theorem \ref{thm:S-Mgen}, we know that there exists
a triple $(C, \bar x, t_{\bar x})$ such that $X_\Omega$ is
isomorphic as a p.p.a.v. to the Jacobian of $C$.

Now 2.b) shows that $U_{\xi_1}=\sigma(U_{\xi_0})$, so the orbit
$\Gamma_V(U_{\xi_0},U_{\xi_1})/\bar \Gamma_V^+$ is finite
dimensional. Theorem~\ref{thm:finitedimNonRam} (see also
Remark~\ref{rem:finiteorbitNRam}) implies that $\sigma$ induces an
automorphism of $C$ satisfying the conditions of 1.
\end{proof}

\begin{rem}
Observe that if the condition 2. of Theorem~\ref{thm:S-Mgen} holds
for one $\xi\in \C^g$, then it holds for every $\xi$ (see
\cite{Shiota}, Theorem 6). Therefore, if the condition 2. of
Theorem~\ref{thm:autram} or \ref{thm:autNram} holds for a given
$\xi_0\in \C^g$, then it holds for every $\xi_0$.
\end{rem}

Finally, we obtain a solution of the Schottky problem for curves
with automorphisms.

\begin{thm}[Characterization]
Let $X_\Omega$ be an irreducible p.p.a.v. of dimension $g>1$.

Then the following conditions are equivalent.
\begin{enumerate}
    \item There exists a projective irreducible smooth curve $C$ of genus $g$ with a non-trivial
    automorphism $\sigma_C:C\to C$ such that $X_\Omega$ is
    isomorphic as a p.p.a.v. to the Jacobian of $C$.

    \item There exist a prime number $p$, $p$ matrices $A^{(1)},\dots, A^{(p)}$ ($A^{(j)}$
    being a $g \times \infty$-matrix of rank $g$) and $p$ symmetric quadratic
    forms $Q^{(1)},\dots, Q^{(p)}$, such that
        \begin{enumerate}
        \item[a)] for some $\xi_0\in \C^g$, the corresponding
        BA-functions satisfy the $(1,\stackrel{p}{\ldots}, 1)$-KP
        hierarchy
            $$
            \res \left( \sum_{j=1}^p z^{-\delta_{ju} - \delta_{jv}} \,
            \psi_{u,\xi_0}^{(j)} (z,t) \psi_{v,\xi_0}^{*(j)} (z,s) \right) \, dz
            \,=\,0
            $$
        \item[b)] there exist $\xi_1\in \C^g$ (depending on $\xi_0$) such that
            $$
            \res \left( \sum_{j=1}^p z^{-\delta_{ju} - \delta_{jv}} \,
            \psi_{v+1,\xi_0}^{(j+1)} (z,\sigma^*(t))
            \psi_{u,\xi_1}^{*(j)} (z,s)\right) \, dz
            \,=\,0
            $$
            where $\sigma^*(t):=
            (t^{(p)}, t^{(1)}, t^{(2)}, \ldots , t^{(p-1)})$.
        \end{enumerate}
\end{enumerate}
\end{thm}

\begin{proof}
Observe that any curve with non-trivial automorphism group admits
an automorphism of prime order $p$ with an orbit consisting of $p$
pairwise distinct points. By the previous theorem, we conclude.

\end{proof}

\begin{rem}
Recall that standard arguments allow us to express the above
equations as an infinite system of partial differential equations
for the $\tau$-function.
\end{rem}



\vskip2truecm
\begin{thebibliography}{MMMM}

\bibitem[ASvM]{ASvM} M. Adler, T. Shiota, P. van Moerbeke, ``Pfaff $\tau$-functions'',
 Math. Ann. {\bf 322} (2002), no. 3, pp. 423 - 476


\bibitem[DJKM]{DJKM} Date, E.; Jimbo, M.; Kashiwara, M. and Miwa, T.,
``Transformation groups for soliton equations'',  Proc. RIMS
Sympos. on Nonlinear Integral Systems, World Scientific, Singapore
(1983), pp. 39-119


\bibitem[F]{Fay}Fay, J.D., ``Bilinear Identities for Theta Functions'',
Preprint

\bibitem[GMP]{GMP} G\'omez Gonz\'alez, E.; Mu\~{n}oz Porras, J. M.; Plaza
Mart\'{\i}n, F. J., ``Prym varieties, curves with automorphisms
and the Sato Grassmannnian'',  Math. Ann. \textbf{327} (2003), no.
4, pp. 609--639

\bibitem[K]{Krichever} Krichever, I.M., ``Methods of algebraic geometry in the
theory of non-linear equations'', Russian Math. Surveys 32:6
(1977), pp 185-213

\bibitem[LM]{LiMulasePrym} Li, Y.; Mulase, M., ``Category of
morphisms of algebraic curves and a characterization of Prym
varieties'', MPI preprint and {\verb"alg-geom/9203002"}

\bibitem[M]{Mulase} Mulase, M., ``Cohomological structure in
soliton equations and Jacobian varieties'', J. Differential Geom.
\textbf{19} (1984), pp. 403--430.

\bibitem[M2]{Mul2} Mulase, M., ``A correspondence between an Infinite
Grassmannian and arbitrary vector bundles on algebraic curves'',
Proc. Symp. Pure Math. \textbf{49} Series A (1989), pp. 39--50.

\bibitem[MP]{MP} Mu\~noz Porras, J.M.; Plaza Mart\'{\i}n, F.J.,
``Equations of the moduli space of pointed curves in the infinite
Grassmannian'', J. Differ. Geom. \textbf{51}  (1999), pp.
431--469.

\bibitem[MP2]{Hurwitz} Mu\~{n}oz Porras, J. M.; Plaza Mart\'{\i}n, F. J.,
``Equations of Hurwitz schemes in the infinite Grassmannian'',
preprint math.AG/0207091.

\bibitem[SS]{Sato} Sato, M.; Sato, Y., ``Soliton equations as dynamical
systems on infinite Grassman manifold'', Lecture Notes Numer.
Appl. Anal. \textbf{5} (1982), pp 259--271.

\bibitem[Sh]{Shiota} Shiota, T., ``Characterization of Jacobian varieties
in terms of soliton equations'', Invent. Math. {\bf 83}
(1986), pp. 333--382.

\bibitem[Sh2]{Sh2} Shiota, T., ``Prym varieties and soliton equations'',
Infinite-dimensional Lie algebras and groups,  Adv. Ser. Math.
Phys., 7, (Luminy-Marseille, 1988), pp. 407--448, World Sci.
Publishing.

\bibitem[SW]{SW} Segal, G.; Wilson, G., ``Loop groups and equations
of KdV type'',  Publ. Math. I.H.E.S. \textbf{61} (1985), pp.
5--64.

\end{thebibliography}
\end{document}